\documentclass[12pt]{article}
\usepackage[T2A]{fontenc}
\usepackage[utf8]{inputenc}
\usepackage[english,russian,ukrainian]{babel}
\usepackage{latexsym,amssymb,amsmath,amsfonts}

\begin{document}
	
	\newcounter{lemma}
	\newcommand{\lemma}{\par \refstepcounter{lemma}%
		{\bf Лема \arabic{lemma}.}}
	
	\newcounter{theorem}
	\newcommand{\theorem}{\par \refstepcounter{theorem}%
		{\bf Теорема \arabic{theorem}.}}
	
	\newcounter{remark}
	\newcommand{\remark}{\par \refstepcounter{remark}%
		{\bf Зауваження \arabic{remark}.}}
	
	\newcounter{corollary}
	\newcommand{\corollary}{\par \refstepcounter{corollary}%
		{\bf Наслідок \arabic{corollary}.}}
	
	\newcounter{definition}
	\newcommand{\definition}{\par \refstepcounter{definition}%
		{\bf Означення \arabic{definition}.}}
	%%%%%%%%%%%%%%%%%%%%%%%%%%%%%%%%%%%%%%%%%%%%%

	%\noindent УДК 517.5  %4, 517.95
	
	\medskip
	
	\vskip 2mm
	%%%%%%%%%%%%%%%%%%%%%%%%%%%%%%%%%%%%%%

\noindent УДК 517.54, 517.95

\medskip

\medskip
{\bf С.В.~Грищук, C.А.~Плакса} (Інститут математики НАН України)

\medskip
{\bf S.V.~Gryshchuk, S.A.~Plaksa} (Institute of Mathematics of NAS of Ukraine, Kyiv,\\ gryshchuk@imath.kiev.ua, plaksa@imath.kiev.ua)
\medskip

\medskip
{\bf Гіперкомплексний метод розв'язання кусково-неперервної бігармонічної задачі в областях з кутовими точками}

\medskip
{\bf A hypercomplex method for solving piecewise continuous biharmonic problem in domains with corner points}

\bigskip
Розглядається кусково-неперервна бігармонічна задача в областях з кутовими точками і відповідна їй крайова задача типу задачі Шварца для моногенних функцій у комутативній бігармонічній алгебрі. Розвинуто метод зведення задач до системи інтегральних рівнянь.

\medskip   A piecewise continuous biharmonic problem in domains with corner points and a corresponding Schwarz type boundary value problem for monogenic functions in a commutative biharmonic algebra are considered. A method for reducing the problems to a system of integral equations is developed.

%%%%%%%%%%%%%%%%%%%%%%%%%%%%%%%%%%%%%%%%%%%%%

\newpage

%	\section{Вступ}
{\bf 1. Вступ.}
	Нехай $D$ --- обмежена область декартової площини $xOy$ з кусково-гладкою
	межею $\partial D$,  яка має скінченний набір
	$\Upsilon:=\left\{\left(x_1, y_1\right), \left(x_2, y_2\right), \dots, \left(x_m,  y_m\right) \right\}$
	кутових точок.
	
	Розглянемо {\em бігармонічну задачу}, яка полягає у відшуканні бігармонічної функції  $u \colon D \longrightarrow \mathbb{R}$,
	що задовольняє бігармонічне рівняння
	\begin{equation}\label{big-eq}
		\Delta^{2}u(x,y):= \left(\frac{\partial^4}{\partial
			x^4}+2\,\frac{\partial^4}{\partial x^2\partial
			y^2}+\frac{\partial^4}{\partial y^4}\right)u(x,y)=0 \qquad \forall (x,y)\in D,
	\end{equation}
	якщо граничні значення її частинних похідних першого порядку
	у точках межі $\partial D$ області $D$, за винятком точок множини $\Upsilon$,
	задовольняють умови
	\begin{equation}\label{big_pr-mCornPoint}
		\begin{array}{l}
			\displaystyle \lim\limits_{(x,y)\to(x_0,y_0),\, (x,y)\in D}
			\frac{\partial
				u(x,y)}{\partial x}=u_{1}(x_0,y_0),\\[4mm]
			\displaystyle \lim\limits_{(x,y)\to(x_0,y_0),\, (x,y)\in D}
			\frac{\partial u(x,y)}{\partial y}=u_{3}(x_0,y_0)\qquad \forall\,
			(x_0,y_0) \in\partial D \setminus \Upsilon,
		\end{array}
	\end{equation}
	і для кожного $j \in \left\{1,2,\dots,m\right\}$  виконуються оцінки
	\begin{equation}\label{ux-nearCornPoints}
		\left|\frac{\partial u(x,y)}{\partial x}\right|+\left|\frac{\partial u(x,y)}{\partial y}\right| \le c\,\left(\left(x-x_j\right)^2+\left(y-y_j\right)^2\right)^{-\alpha_j /2},
		\quad (x,y)\to (x_j,y_j)\,, \,  (x,y)\in D\setminus \Upsilon\,,\\
	\end{equation}
	де $\alpha_j\in (0,1)$, $j=\overline{1,m}$, і стала\, $c$\, не залежить від $(x,y)$.
	
	Будемо припускати, що задані функції $u_j \colon \partial D\setminus\Upsilon\longrightarrow\mathbb{R}$, $j\in\{1,3\}$, --- неперервні на множині $\partial D\setminus\Upsilon$. У той же час, з оцінок \eqref{ux-nearCornPoints} випливає, що для них при $j=\overline{1,m}$ мають виконуватися оцінки
	\begin{equation} \label{u_j-oz-nearCornPoints}
		\left|u_j(x,y)\right| \le c\,\left(\left(x-x_j\right)^2+\left(y-y_j\right)^2\right)^{-\alpha_j /2},
		\quad (x,y)\to (x_j,y_j)\,, \,  (x,y)\in\partial D\setminus\Upsilon\,,
	\end{equation}
	де стала\, $c$\, не залежить від $(x,y)$.

	Класична постановка бігармонічної задачі (див., наприклад, \cite{Mikhlin}) передбачає, що
	$\mathcal{X}=\emptyset$, а межа  $\partial D$  % $\partial D_z$ 
	є гладкою кривою. Тоді граничні умови  \eqref{big_pr-mCornPoint}  розглядаються
	для усіх $\left(x_0,y_0\right)\in \partial D$, а оцінки \eqref{ux-nearCornPoints} та \eqref{u_j-oz-nearCornPoints} відсутні.
	
	Теорія крайових задач для еліптичних рівнянь в областях з кусково-гладкими межами розвинута в роботах багатьох авторів (див., наприклад, роботи \cite{Kondrat'ev-TMMO-67, Kufner-Sandig,MaNazPlam-ASumptTh-V1} та літературу в них), при цьому вивчено питання існування і єдиності розв'язків крайових задач для еліптичних операторів у вагових просторах, досліджено асимптотичні властивості розв'язків в околах кутових точок та вплив кутових точок на гладкість розв'язків тощо.
	
	Для розв'язання крайових задач для бігармонічних функцій
	розвинені методи, що використовують аналітичні функції комплексної змінної, техніка використання яких базується на зображенні бігармонічних функцій формулою Гурса. Це дозволяє зводитити крайові задачі для бігармонічних функцій до відповідних крайових задач для пари аналітичних функцій. Далі з використанням зображення аналітичних функцій інтегралами типу Коші в загальному випадку отримують систему інтегро-диференціальних рівнянь. У випадку, коли межа області є кривою Ляпунова, зазначена система, як правило, зводиться до системи рівнянь Фредгольма. Така схема розроблена (див., наприклад, \cite{Mush_upr,Lurie_engl,Mikh_int_eq_Th-El}) для розв'язання основних задач плоскої
	теорії пружності з використанням спеціальної бігармонічної функції, яка називається функцією напружень Ері.
	Інші методи зведення крайових задач плоскої теорії пружності до інтегральних рівнянь розвинені в роботах
	\cite{Mikhlin,Mikh_kniga,Kupragze63,Lopatinskii,Maz'ya91Engl}.
	
	У роботі \cite{Magnaradze-38} показано, що для обмеженої області, межа якої є кривою з обмеженим обертанням, тобто кривою Радона  (див. \cite{Radon-46}), і з кутовими точками, у випадку заданих граничних функцій, що задовольняють умову Гельдера, основні крайові задачі плоскої теорії пружності можна звести до інтегрального рівняння типу рівняння Мусхелішвілі (див. \cite{Mush_upr}).
	Розглядалися також крайові задачі для бігармонічних функцій в областях з конкретними кусково-гладкими межами (див., наприклад, \cite{Polozhyi,G-Albinus-83PAN,Meln-PLaksa-PBP-quadrante,GrPl-Angle-umz}).
	
	У роботах \cite{IJPAM_13,mon-f-bih-BVP-MMAS,Gr-1Fr17,BeghGrPl19,GrPlHBVMbF-Flaut21,SrPlCanD-UMB21}
	розивинуто метод розв'язання бігармонічної задачі, який базується на зв'язку бігармонічних функцій з моногенними функціями у
	двовимірній комутативній асоціативній алгебрі над полем комплексних чисел та
	на пред\-став\-лен\-ні розв'язків гіперкомплексними інтегралами, ана\-ло\-гіч\-ни\-ми до класичних
	інтегралів Шварца і інтегралів типу Коші.
	Цей гіперкомплексний метод у загальному випадку дозволяє редукувати бігармонічну задачу
	безпосередньо до системи інтегральних рівнянь, оминаючи інтегро-диференціальні рівняння.
	У роботах \cite{mon-f-bih-BVP-MMAS,Gr-1Fr17,BeghGrPl19,GrPlHBVMbF-Flaut21},
	встановлено достатні умови фредгольмовості вказаної системи для обмежених областей з гладкими межами, що належать більш широким класам,
	ніж клас кривих Ляпунова, який, як правило, розглядався раніше в
	плоскій теорії пружності.
	
	У даній роботі гіперкомплексний метод розвинуто для зведення кусково-неперервної бігармонічної задачі для обмеженої області з кутовими точками до системи інтегральних рівнянь.

%	\section{Крайова задача для моногенних функцій, асоційована з бігармонічною задачею для обмеженої області з кутовими точками}
{\bf 2. Крайова задача для моногенних функцій, асоційована з бігармонічною задачею для обмеженої області з кутовими точками.}
		У роботі \cite{KM-BFf} асоціативну комутативну двовимірну алгебру $\mathbb{B}$ з одиницею $e$ над полем комплексних чисел $\mathbb{C}$
	названо {\em бігар\-моніч\-ною},
	якщо вона містить базис $\{e_1,e_2\}$, що задовольняє вимоги
	\begin{equation}\label{biharm-bas}
		(e_1^2+e_2^2)^2=0,\qquad e_1^2+e_2^2\ne 0
	\end{equation}
	(який також названо {\em бігармонічним}), і запропоновано наступну таблицю множення для такого базису:
	\begin{equation}\label{tabl-umn}
		e_1^2=e_1,\quad e_2e_1=e_2,\quad e_2^2=e_1+2ie_2,
	\end{equation}
	де $i$~--- уявна комплексна одиниця.
	
	У роботі  \cite{Mel86} доведено єдиність бігармонічної алгебри $\mathbb{B}$ і показано, що вона породжується  небігармонічним базисом  $\{e_1,\rho\}$, де
	\begin{equation} \label{rho}
		\rho=2e_1+2ie_2\,,
	\end{equation}
	при цьому\, $\rho^{2}=0$\,,
	а також описано усі бігармонічні базиси в $\mathbb{B}$.
	Зауважимо, що алгебра $\mathbb B$ ізоморфна
	чотиривимірним  алгебрам над полем дійсних чисел
	$\mathbb{R}$, розглянутим у роботах \cite{Sodbero,Douglis-53}.
	
	Задамо евклідову норму $\|a\|:=\sqrt{|z_1|^2+|z_2|^2}$ в алгебрі  $\mathbb{B}$, де\, $a=z_1e_1+z_2e_2$ і  $z_1, z_2\in \mathbb{C}$.
	
	Як і в роботі \cite{KM-BFf}, розглянемо  {\it бігармонічну площину} $\mu_{e_1,e_2}:=\{\zeta=x\,e_1+y\,e_2 :
	x,y\in\mathbb R\}$.
	
	Області $D$ декартової площини $xOy$ поставимо у відповідність конгруентну їй
	область $D_{\zeta}:=\{\zeta=xe_1+ye_2 : (x,y)\in D\}$ у площині
	$\mu_{e_1,e_2}$, а її замиканню $\overline{D}$~--- замикання $\overline{D_{\zeta}}$.
	Їх межі позначатимемо відповідно через $\partial D$ та $\partial D_{\zeta}$.
	
	Скрізь надалі\, $\zeta:=x\,e_1+y\,e_2$, $z:=x+iy$,
	% \equiv \mathrm{Re} \, z + i\, \mathrm{Im} \, z$, 
	де $(x,y)\in D$, і\,
	$\zeta_{0}:=x_{0}e_{1} + y_{0}e_{2}
	$, де  $(x_0,y_0)\in \partial D$.
	
	%$Z_0=x_0 + i y_0$, де
	
	Оскільки в бігармонічній площині відсутні дільники нуля, то
	похідна функції $\Phi \colon D_{\zeta}\longrightarrow \mathbb{B}$ визначається так, як і для аналітичних функцій комплексної змінної, а саме:
	\[\Phi'(\zeta):=\lim\limits_{h\to 0,\, h\in\mu}
	\bigl(\Phi(\zeta+h)-\Phi(\zeta)\bigr)\,h^{-1}\,.\]
	Функція  $\Phi \colon D_{\zeta}\longrightarrow \mathbb{B}$ називається
	{\em моногенною} в області  $D_{\zeta}$, якщо її похідна
	$\Phi'(\zeta)$ існує у кожній точці  $\zeta\in D_{\zeta}$.

	Кожна функція $\Phi\colon D_{\zeta}\longrightarrow \mathbb{B}$ подається у вигляді
	\begin{equation}\label{mon-funk}
		\Phi(\zeta)=
		U_{1}(x,y)\,e_1+U_{2}(x,y)\,ie_1+U_{3}(x,y)\,e_2+U_{4}(x,y)\,ie_2\,,
	\end{equation}
	де  $U_{l}\colon D\longrightarrow \mathbb{R}$,  $l=\overline{1,4}$,
	--- дійснозначні компоненти-функції.
	Для них також будемо  використовувати позначення
	$\mathrm{U}_{l}\left[\Phi\right]:=U_{l}$, $l=\overline{1,4}$.

	У роботі \cite{KM-BFf} доведено, що функція   $\Phi\colon
	D_{\zeta}\longrightarrow \mathbb{B}$ --- моногенна  в області $D_{\zeta}$ тоді і тільки тоді, коли всі її дійснозначні компоненти з розкладу
	(\ref{mon-funk}) диференційовні в  $D$ і виконується наступний аналог умов Коші--Рімана:
	\begin{equation}\label{usl_K_R}
		\frac{\partial \Phi(\zeta)}{\partial y}=\frac{\partial
			\Phi(\zeta)}{\partial x}\,e_2.
	\end{equation}

	У роботах \cite{GrPl_umz-09,dopovidi_09,Conrem_13} доведено, що кожна моногенна функція $\Phi \colon D_{\zeta}\longrightarrow
	\mathbb{B}$ має похідні  $\Phi^{(n)}(\zeta)$ усіх порядків
	в області  $D_{\zeta}$ і задовольняє рівняння  \eqref{big-eq}
	в силу першого зі співвідношень  (\ref{biharm-bas}) і рівності
	\[\Delta^{2}\Phi(\zeta)=\Phi^{(4)}(\zeta)\,(e_1^2+e_2^2)^2.\]
	Тому всі компоненти  $U_{l}\colon D\longrightarrow \mathbb{R}$,
	$l=\overline{1,4}$, з розкладу  \eqref{mon-funk} є бігармонічними функціями в області  $D$.

	Крім того, кожна бігармонічна в  $D$
	функція  $U(x,y)$ є першою компонентою  $U_{1}\equiv U$ з розкладу  (\ref{mon-funk}) для деякої
	моногенної функції
	$\Phi \colon D_{\zeta}\longrightarrow \mathbb{B}$. Усі такі  функції  $\Phi$
	знайдено в роботах \cite{GrPl_umz-09,Conrem_13} у явному вигляді для обмежених однозв'язних областей.
	
	Нехай $\Phi_{1}$ ---  моногенна функція в області $D_{\zeta}$, яка має своєю першою компонентою
	шукану функцію $u(x,y)$ бігармонічної задачі, тобто $\Phi_{1}$ має розклад
	\[\Phi_{1}(\zeta)=u(x,y)\,e_1+U_{2}(x,y)\,ie_1+U_{3}(x,y)\,e_2+U_{4}(x,y)\,ie_2 \qquad \forall\,\zeta \in D_{\zeta}\,.\]
	
	З умови (\ref{usl_K_R}) для $\Phi=\Phi_1$ випливає рівність $\partial U_3(x,y)/\partial x=\partial u(x,y)/\partial y$. Тому
	\begin{equation*}
		\Phi_{1}'(\zeta)=\frac{\partial u (x,y)}{\partial x}\,e_1+
		\frac{\partial U_{2}(x,y)}{\partial x}\,ie_1+ \frac{\partial
			u(x,y)}{\partial y}\,e_2+ \frac{\partial U_{4}(x,y)}{\partial
			x}\,ie_2  \qquad \forall\,\zeta \in D_{\zeta}\,.
	\end{equation*}
	
	Отже, бігармонічна задача з крайовими умовами \eqref{big_pr-mCornPoint} редукується до еквівалентної крайової задачі про знаходження моногенної в області
	$D_{\zeta}$ функції $\Phi\equiv\Phi_{1}'$ у випадку, коли значення двох компоненент $\mathrm{U}_{l}\left[\Phi\right]$, $l\in
	\{1,3\}$, розкладу (\ref{mon-funk}) задані на множині $\partial D_{\zeta}\setminus \Upsilon_{\zeta}$,
	де $\Upsilon_{\zeta}:=\left\{\zeta_1, \zeta_2, \dots, \zeta_m \right\}$ ---
	набір кутових точок $\zeta_j:= e_1 \, x_j +  e_2\, y_j$, $j=\overline{1,m}$,  межі $\partial D_{\zeta}$.
	
	Дамо точне формулювання крайової задачі для моногенних функцій, що є предметом нашого розгляду.
	
	Розглянемо кусково-неперервну (1-3)-{\it задачу} про знаходження моногенної функції $\Phi \colon
	D_{\zeta}\longrightarrow \mathbb{B}$, для якої граничні значення компонент $\mathrm{U}_{l}
	\left[\Phi(\zeta)\right]$, $l\in \{1,3\}$, з розкладу (\ref{mon-funk}) задовольняють крайові умови
	\begin{equation}\label{Pr13week}
		\lim\limits_{\zeta\to\zeta_0,\zeta\in D_{\zeta}} \mathrm{U}_{l}
		\left[\Phi(\zeta)\right]=u_{l}(\zeta_0) \qquad \forall\,\zeta_0\in \partial D_{\zeta}\setminus\Upsilon_{\zeta},\qquad
		l\in\{1,3\},
	\end{equation}
	де $u_{l} \colon \partial D_{\zeta} \setminus \Upsilon_{\zeta} \longrightarrow \mathbb{R}$, $l\in \{1,3\}$, --- задані неперервні на
	$\partial D_{\zeta}\setminus\Upsilon_{\zeta}$ функції, які  ототожнені з відповідними граничними функціями
	бігармонічної задачі:
	$u_{l}(\zeta_0)\equiv u_l(x_0,y_0)$, $l\in\{1,3\}$.
	Крім того, шукана функція $\Phi$  у відповідності з оцінками \eqref{ux-nearCornPoints} має задо\-воль\-ня\-ти для кожних $l \in\{1,3\}$ і $j \in \left\{1,2,\dots,m\right\}$  оцінки
	\begin{equation}\label{Ul-toCornPoin}
		\Bigl|\mathrm{U}_{l}\left[\Phi(\zeta)\right]\Bigr|  \le c\,\|\zeta-\zeta_j \|^{-\alpha_j},\qquad  \zeta\to \zeta_j, \quad  \zeta\in
		D_{\zeta},
	\end{equation}
	де $\alpha_j \in (0,1)$ і стала\, $c$\, не залежить від $\zeta$.
	
	У роботі \cite{Kovalov}  задачі такого типу названо {\it бігармонічними задачами Шварца} з огляду на їх пряму аналогію з класичною задачею Шварца про відшукання аналітичної функції комплексної змінної, коли  значення її дійсної частини задано на межі області.
	
	Після розв'язання  (1-3)-задачі моногенна функція $\Phi_{1}$ (а разом з нею і перша її компонента --- розв'язок  $u$ бігармонічної  задачі)
	знаходиться в результаті контурного інтегрування
	\[\Phi_{1}(\zeta)=\int\limits_{\gamma_{\zeta_1, \zeta}}\Phi(\tau)\,d\tau+{\rm const} \qquad \forall\,\zeta\in D_{\zeta} \]
	взовж спрямлюваної кривої $\gamma_{\zeta_1, \zeta}\subset D_{\zeta}$, що сполучає фіксовану точку $\zeta_1\in D_{\zeta}$ з точкою $\zeta$.
	Зазначимо, що в силу справедливості інтегральної теореми Коші для моногенних функцій (див. \cite{Conrem_13})
	результат такого інтегрування не залежить від вибору кривої $\gamma_{\zeta_1, \zeta}$.

	Будемо шукати розв'язки  (1-3)-задачі у класі функцій, що подаються у вигляді бігармонічного інтеграла типу Коші
	\begin{equation}\label{cl-sol-cone}
		\Phi(\zeta)= \frac{1}{2\pi i}\int\limits_{\partial
			D_\zeta}\varphi(\tau)(\tau-\zeta)^{-1}\,d\tau
		%=:\mathcal{B}[\varphi](\zeta)
\qquad \forall\, \zeta \in D_{\zeta}\,,
	\end{equation}
	де
	\begin{equation}\label{varphi13}
		\varphi(\zeta)=\varphi_{1}(\zeta)\,e_1+\varphi_{3}(\zeta)\,e_2 \quad \forall\, \zeta \in \partial D_{\zeta},
	\end{equation}
	і функції
	$\varphi_{l} \colon
	\partial D_{\zeta}\setminus\Upsilon_{\zeta} \longrightarrow \mathbb{R}$, $l\in\{1,3\}$, --- неперервні на
	$\partial D_{\zeta}\setminus\Upsilon_{\zeta}$,
	а умови на їх поведінку в околі точок набору $\Upsilon_{\zeta}$ будуть наведені нижче.

%	\section{Допоміжні твердження}
{\bf 3. Допоміжні твердження.}
	
%	\subsection{Кусково-неперервне продовження на межу інтеграла типу Коші, щільність якого є функцією двох комплексних змінних}
{\bf 3.1. Кусково-неперервне продовження на межу інтеграла типу Коші, щільність якого є функцією двох комплексних змінних.}	
	Нехай $D_{z}:=\{z=x+iy : \left(x,y\right)\in D\}$ --- область комплексної площини  $\mathbb{C}$, конгруентна області $D$,
	і $z_j:=x_j+iy_j$, $j=\overline{1,m}$, --- кутові точки межі $\partial D_z$ області $D$.
	
	Розглянемо конформне відображення  $z= \sigma(Z)$ одиничного круга
	$\mathcal{U}:= \{Z\in\mathbb{C}: |Z|<1\}$ на область  $D_{z}$,
	яке очевидно визначає також гомеоморфізм одиничного кола $\Gamma:=\{S\in\mathbb{C}: |S| = 1\}$ та межі $\partial D_{z}$.
	
	Нехай $\mathcal{X}:= \{X_1,X_2, \dots, X_m\}$ --- набір точок кола $\Gamma$ таких, що $\sigma(X_j)=z_j$, $j=\overline{1,m}$.
	
	Позначимо $\mathcal{R}_{0} 	:= \min\limits_{ X_{k}, X_{j} \in \mathcal{X},\,
		k  \ne j}\left\{\left|X_k-X_{j}\right|\right\}$, якщо $m>1$,  і  $\mathcal{R}_{0}:=1$ у випадку $m=1$.
	
	Для кожної точки $Z$, що належить замиканню $\overline{\mathcal{U}}$ круга $\mathcal{U}$, позначимо
	\[ \mathcal{R}(Z):= \min\limits_{X_j \in \mathcal{X}}\left\{\left|Z-X_{j}\right|\right\}.\]

	Для $r>0$, $E\subset\mathbb{C}$ і $T_0\in \mathbb{C}$ позначимо
	$E_{r}\left(T_0\right):=\left\{T\in E: \left|T-T_0\right|\le r\right\}$,
	а також
	\[ E_{r}\left\langle \mathcal{X}\right\rangle := \bigcup\limits_{j=1}^{m} E_{r} \left(X_j \right). \]

	Сингулярний інтеграл розуміємо у сенсі головного значення за Коші
	\[\int\limits_{\Gamma}\frac{\Omega(S,\,\cdot\,)}{S-T}\,dS:=\lim\limits_{\delta\to 0+0}\,\int\limits_{\Gamma\setminus\Gamma_{\delta}(T)}
	\frac{\Omega(S,\,\cdot\,)}{S-T}\,dS, \qquad T\in\Gamma, \]
	за умови, що границя в правій частині рівності існує.

	\medskip
	\begin{lemma}\label{contBikhIntCorPs}
		{\it Нехай  функція $\Omega \colon \left(\Gamma \setminus \mathcal{X}\right) \times \left(\overline{\mathcal{U}} \setminus \mathcal{X} \right) \longrightarrow \mathbb{C} $ задовольняє оцінки
			\begin{equation}
				\label{omeg-loc-bound}
				\left| \Omega\left(S,Z\right)\right| \le \frac{c}{\left(\mathcal{R}(S)\right)^{\beta}} \,
				\min\left\{\frac{\omega_{0}\left(|S-Z|\right)}
				{\omega_{0}\left(\frac{1}{2}\mathcal{R}(Z)\right)}, \,  1 \right\} \qquad
				\forall Z \in \mathcal{U} \quad \forall S \in \Gamma \setminus \mathcal{X},
			\end{equation}
			\begin{equation}
				\label{omeg-locBound-circle}
				\left| \Omega\left(S,T\right)\right| \le \frac{c}{\left(\mathcal{R}(S)\right)^{\beta}} \,
				\min\left\{\frac{\omega_{1}\left(|S-T|\right)}
				{\omega_{1}\left(\frac{1}{2}\mathcal{R}(T)\right)}, \,  1 \right\} \qquad
				\forall S, T \in \Gamma \setminus \mathcal{X},
			\end{equation}
			\begin{multline} \label{omeg-locBound-diff}
				\left| \Omega\left(S,Z\right)- \Omega\left(S,Z_0 \right) \right| %\\
				\le \frac{c  }{\left(\mathcal{R}(S)\right)^{\beta} \left(\mathcal{R}(Z_0)\right)^{\beta_0}}
				\min\left\{\frac{\omega_{2}\left(|S-Z_0|\right)}
				{\omega_{2}\left(\frac{1}{2}\mathcal{R}(Z_0)\right)}, \, 1 \right\}
				\frac{|Z-Z_0|}{|S-Z_0|}\\[2mm]
				\forall Z_0 \in \Gamma \setminus \mathcal{X} \quad \forall Z \in \overline{\mathcal{U}}_{\frac{1}{4}\mathcal{R}(Z_0)}(Z_0) \setminus \mathcal{X} \quad
				\forall S \in \Gamma \setminus \mathcal{X}: |S-Z_0| \ge 2|Z-Z_0|,
			\end{multline}
			де стала $c$ не залежить від $S$, $Z$, $T$, $Z_0$\emph{;} $\beta \in (0; 1)\emph{;} \, \beta_0 \ge 0$\emph{;}
			$\omega_k \colon (0; +\infty) \longrightarrow (0; +\infty)$ при $k=0,1,2$ --- неспадні обмежені  функції, які задовольняють умови\emph{:}
			\[\omega_k(\varepsilon)\to 0, \quad \varepsilon \to 0, \qquad k\in\{0,2\},\]
			%функція $\omega_1$ задовольняє умову Діні
			\begin{equation}\label{dini}
				\int\limits_{0}^{2}\frac{\omega_1(\eta)}{\eta}\, d\eta < \infty.
			\end{equation}
			
			Тоді 
			\begin{equation}\label{limChOmeg-CP}
				\lim\limits_{Z\to Z_{0},\, Z \in \mathcal{U}}\int\limits_{\Gamma}\frac{\Omega\left(S,Z\right)}{S-Z}\, dS =  \int\limits_{\Gamma}\frac{\Omega(S,Z_{0})}{S-Z_{0}}\, dS \qquad \forall Z_{0}\in \Gamma \setminus \mathcal{X}.
			\end{equation}
		}
	\end{lemma}
	
	{\bf \em Доведення.}
	Нехай $Z_{0}\in \Gamma \setminus \mathcal{X}$, $Z \in \mathcal{U}$ і $\varepsilon:= \left|Z-Z_{0} \right| \le \frac{1}{8}\, \mathcal{R}\left(Z_0\right)$. Введемо в розгляд точку
	$Z_{\ast}\in \Gamma$  таку, що
	$\left|Z_{\ast}-Z\right|= \min\limits_{S \in \Gamma}\left|S -Z\right|$.
	
	Розглянемо різницю
	\[ \int\limits_{\Gamma}
	\frac{\Omega(S,Z)}{S-Z}\, dS-  \int\limits_{\Gamma}
	\frac{\Omega(S,Z_{0})}{S-Z_{0}}\, dS=
	\int\limits_
	{\Gamma_{
			\frac{1}{2} \mathcal{R}(Z_0)}
		\left\langle \mathcal{X}\right\rangle}
	\frac{\Omega\left(S,Z\right)- \Omega\left(S,Z_{0}\right)}
	{S-Z} \, d S +\]
	\[+  \int\limits_{\Gamma_{\frac{1}{2} \mathcal{R}(Z_0)}\left\langle \mathcal{X}\right\rangle}
	\Omega(S,Z_{0})\left(\frac{1}{S-Z}- \frac{1}{S- Z_{0}}\right)\, dS+\int\limits_{\Gamma  \setminus\,
		\Gamma_{\frac{1}{2} \mathcal{R}(Z_0)}\left\langle \mathcal{X}\right\rangle}
	\frac{\Omega\left(S,Z\right)- \Omega\left(S,Z_{\ast}\right)}{S-Z}\,d S +\]
	\[ +\int\limits_{\Gamma  \setminus\,
		\Gamma_{\frac{1}{2} \mathcal{R}(Z_0)}\left\langle \mathcal{X}\right\rangle}
	\frac{\Omega\left(S,Z_{\ast}\right)}{S-Z}\, dS
	- \int\limits_{\Gamma  \setminus\,
		\Gamma_{\frac{1}{2} \mathcal{R}(Z_0)}\left\langle \mathcal{X}\right\rangle}
	\frac{\Omega\left(S,Z_{0}\right)}{S-Z_0}\, dS
	=\]
	\[	=:  j_{1} + j_{2} + j_{3} + j_{4}-j_{5}.\]
		
Далі у доведенні через\, $c$\, позначено сталі, значення яких не залежить  від $S$, $Z$, $Z_*$ і $Z_0$,
	але, взагалі кажучи, різні навіть у межах одного ланцюжка нерівностей.
	
	Враховуючи співвідношення
	\begin{equation}\label{estim-IntGamZ01Int}
		|S-Z|\ge \left({3}/{8}\right)\mathcal{R}(Z_0), \quad |S-Z_0| \ge \left(1/2\right) \mathcal{R}(Z_0)
		\qquad \forall   S\in\Gamma_{\frac{1}{2}\mathcal{R}\left(Z_0\right)} \left\langle \mathcal{X}\right\rangle
	\end{equation}
	і нерівність
	\eqref{omeg-locBound-diff},  одержуємо оцінку інтеграла $j_{1}$:
	\[\left|j_{1}\right| \le \int\limits_
	{\Gamma_{\frac{1}{2} \mathcal{R}(Z_0)} \left\langle \mathcal{X}\right\rangle}
	\left|
	\frac{\Omega\left(S,Z\right)- \Omega\left(S,Z_{0}\right)}{S-Z}\right|\, |d S| \le \]
	\[ \le
	c \, \frac{8}{3 \, \mathcal{R}(Z_0)} \int\limits_
	{\Gamma_{\frac{1}{2} \mathcal{R}(Z_0)}\left\langle \mathcal{X}\right\rangle}
	\frac{\varepsilon\, \min\left\{\frac{\omega_{2}\left(|S-Z_0|\right)}
		{\omega_{2}\left(({1}/{2})\,\mathcal{R}(Z_0)\right)}, \, 1 \right\}}{\left(\mathcal{R}(S)\right)^{\beta} \left(\mathcal{R}(Z_0)\right)^{\beta_0}\, |S-Z_0|} |d S| \le \]
	\[ \le c \frac{\varepsilon}
	{\left(\mathcal{R}(Z_0)\right)^{2+\beta_0}} \sum_{j=1}^{m} \, \int\limits_
	{\Gamma_{\frac{1}{2} \mathcal{R}(Z_0)\left(X_j\right)}}  \frac{|dS|}{\left(\mathcal{R}(S)\right)^{\beta}} \le
	c \frac{\varepsilon}
	{\left(\mathcal{R}(Z_0)\right)^{2+\beta_0}} \sum_{j=1}^{m} \, \int\limits_
	{0}^{\mathcal{R}(Z_0)/2}  \frac{d\eta }{\eta^{\beta}}\le
	c\frac{\varepsilon}{\left(\mathcal{R}(Z_0)\right)^{1+\beta+\beta_0}}. \]
	
	Враховуючи нерівності
	\eqref{estim-IntGamZ01Int} і \eqref{omeg-locBound-circle}, отримуємо оцінку інтеграла $j_{2}$:
	\[\left|j_{2}\right| \le    c \frac{\varepsilon}{\left(\mathcal{R}(Z_0)\right)^{2}} \,
	\sum_{j=1}^{m} \, \int\limits_
	{\Gamma_{\frac{1}{2} \mathcal{R}(Z_0)\left(X_j\right)}}
	\frac{|d S|}{\left(\mathcal{R}(S)\right)^{\beta}}
	\le
	c \frac{\varepsilon}
	{\left(\mathcal{R}(Z_0)\right)^{1 + \beta}}. \]

	Позначимо $\delta:= \left|Z-Z_{\ast}\right|\equiv 1- |Z|$.
	Тоді маємо
	\[\left|j_{3}\right| \le \int\limits_{\Gamma_{2\delta}\left(Z_{\ast}\right)}
	\frac{\left|\Omega\left(S,Z\right)\right|}{|S-Z|} \left|d S\right| +
	\int\limits_{\Gamma_{2 \delta}\left(Z_{\ast}\right)}
	\frac{\left|\Omega\left(S,Z_{\ast}\right)\right|}{|S-Z|} \left|d S\right|+\]
		\[	+ \int\limits_
	{\left(\Gamma \setminus \Gamma_{2 \delta}\left(Z_{\ast}\right)\right)  \setminus\,
		\Gamma_{\frac{1}{2} \mathcal{R}(Z_0)}\left\langle \mathcal{X}\right\rangle }
	\frac{\left|\Omega\left(S,Z\right) - \Omega\left(S,Z_{\ast}\right)\right|}{|S-Z|} \left|d S\right|=: j_{3, 1}+ j_{3, 2} +j_{3,3}.\]
		
	Враховуючи оцінки \eqref{omeg-loc-bound}, \eqref{omeg-locBound-circle} та співвідношення $\mathcal{R}(Z)\ge (7/8)\mathcal{R}(Z_0)$,\,
	$\mathcal{R}(Z_{\ast})\ge (3/4)\mathcal{R}(Z_0)$,\,
	$\delta \le |S-Z| \le 3\delta \le 3 \varepsilon$,\, $|S-Z_{\ast}|\le 2\delta$ і  $\mathcal{R}(S) \ge (1/2) \mathcal{R}(Z_0)$ при всіх $S\in \Gamma_{2 \delta}\left(Z_{\ast}\right)$,\,
	послідовно отримуємо
	\[ \left|j_{3, 1}\right| \le
	c\int\limits_{\Gamma_{2 \delta}\left(Z_{\ast}\right)}
	\frac{1}{|S-Z|} \frac{\omega_{0}\left(|S-Z|\right)}
	{\omega_{0}\left(\frac{1}{2}\mathcal{R}(Z)\right)}  \frac{\left|d S\right|}{\left(\mathcal{R}(S)\right)^{\beta}} \le \]
	\[\le c\,
	\frac{ \omega_{0}\left(3\delta\right)}{\omega_{0}\left(\frac{7}{16}\mathcal{R}(Z_0)\right) \left(\mathcal{R}\left(Z_0\right)\right)^{\beta}\,\delta}\, \int\limits_{\Gamma_{2 \delta}\left(Z_{\ast}\right)}
	{\left|d S\right|}
	\le c\,
	\frac{ \omega_{0}\left(3\varepsilon\right)}{\left(\mathcal{R}\left(Z_0\right)\right)^{\beta} \omega_{0}\left(\frac{7}{16}\mathcal{R}(Z_0)\right) }, \]
	
	\[ \left|j_{3, 2}\right| \le
	c\int\limits_{\Gamma_{2 \delta}\left(Z_{\ast}\right)}
	\frac{1}{|S-Z|} \frac{\omega_{1}\left(|S-Z_{\ast}|\right)}
	{\omega_{1}\left(\frac{1}{2}\mathcal{R}(Z_{\ast})\right)}  \frac{\left|d S\right|}{\left(\mathcal{R}(S)\right)^{\beta}} \le \]
	\[\le c\,
	\frac{ \omega_{1}\left(2\delta\right)}{\omega_{1}\left(\frac{3}{8}\mathcal{R}(Z_0)\right) \left(\mathcal{R}\left(Z_0\right)\right)^{\beta}\,\delta}\, \int\limits_{\Gamma_{2 \delta}\left(Z_{\ast}\right)}
	{\left|d S\right|}
	\le c\,
	\frac{ \omega_{1}\left(2\varepsilon\right)}{\left(\mathcal{R}\left(Z_0\right)\right)^{\beta} \omega_{1}\left(\frac{3}{8}\mathcal{R}(Z_0)\right) }. \]
	
	Враховуючи оцінку \eqref{omeg-locBound-diff} та співвідношення
	$\mathcal{R}(Z_{\ast})\ge (3/4)\mathcal{R}(Z_0)$, $|Z-Z_{\ast}|=\delta$,
	$\mathcal{R}(S)\ge (1/2)\mathcal{R}(Z_0)$ і $|S-Z|\ge (1/2)\left|S- Z_{\ast}\right|$ при всіх $S\in \left(\Gamma \setminus \Gamma_{2 \delta}\left(Z_{\ast}\right)\right)  \setminus\,
	\Gamma_{\frac{1}{2} \mathcal{R}(Z_0)}\left\langle \mathcal{X}\right\rangle$, отримуємо
	\[ \left|j_{3, 3}\right| \le c\, \frac{\delta}{\omega_{2}\left(\frac{3}{8}\mathcal{R}(Z_0)\right)\left(\mathcal{R}(Z_0)\right)^{\beta+\beta_0}}
	\int\limits_
	{\left(\Gamma \setminus \Gamma_{2 \delta}\left(Z_{\ast}\right)\right)  \setminus\,
		\Gamma_{\frac{1}{2} \mathcal{R}(Z_0)}\left\langle \mathcal{X}\right\rangle }
	\frac{\omega_{2}\left(|S-Z_{\ast}|\right)}{|S-Z_{\ast}|^2}|dS|\le\]
	\[\le  c\, \frac{1 }{\omega_{2}\left(\frac{3}{8}\mathcal{R}(Z_0)\right)\left(\mathcal{R}(Z_0)\right)^{\beta+\beta_0}} \, \delta \int\limits_{2\delta}^{2}\frac{\omega_{2}\left(\eta\right)}{\eta^2} d\eta.\]
	
	Подамо різницю $j_{4}-j_{5}$ у вигляді
	\[j_{4}-j_{5}= \int\limits_{\Gamma_{2 \varepsilon}\left(Z_{0}\right)}
	\frac{\Omega\left(S,Z_{\ast}\right)}{S-Z}\, d S -
	\int\limits_{\Gamma_{2 \varepsilon}\left(Z_{0}\right)}
	\frac{\Omega\left(S,Z_{0}\right)}{S-Z_{0}}\, d S + \]
	\[+ \int\limits_
	{\Gamma \setminus \left(\Gamma_{\frac{1}{2} \mathcal{R}(Z_0)}\left\langle \mathcal{X}\right\rangle \, \bigcup \, \Gamma_{2 \varepsilon}\left(Z_{0}\right)\right)}
	\frac{\Omega\left(S,Z_{\ast}\right) - \Omega\left(S,Z_{0}\right)}{S-Z} \, d S + \]
	\[	+ \left(Z-Z_{0}\right) \int\limits_
	{ \Gamma \setminus \left(\Gamma_{\frac{1}{2} \mathcal{R}(Z_0)}\left\langle \mathcal{X}\right\rangle \, \bigcup \, \Gamma_{2 \varepsilon}\left(Z_{0}\right)\right)}
	\frac{\Omega\left(S,Z_{0}\right)}{\left(S-Z_0\right)\left(S-Z\right)} \, d S=: i_{1}- i_{2}+ i_{3} +i_{4}.\]
		
	Враховуючи
	нерівність \eqref{omeg-locBound-circle} та співвідношення $\mathcal{R}(Z_{\ast})\ge (3/4)\mathcal{R}(Z_0)$,  $\mathcal{R}(S)\ge (3/4)\mathcal{R}(Z_0)$ при всіх $S\in \Gamma_{2 \varepsilon}\left(Z_0\right)$, $|S-Z|\ge (1/2)\left|S- Z_{\ast}\right|$ для кожного $S\in \Gamma$,   отримуємо
	\[\left|i_{1}\right| \le 2 \int\limits_{\Gamma_{2 \varepsilon}\left(Z_{0}\right)}
	\frac{\left|\Omega\left(S,Z_{\ast}\right)\right|}{\left|S-Z_{\ast}\right|}\, \left|d S\right| \le \frac{c}{\left(\mathcal{R}\left(Z_{0}\right)\right)^{\beta} \omega_{1}\left(({1}/{2})\,\mathcal{R}\left(Z_{\ast}\right)\right)}
	\int\limits_{\Gamma_{2 \varepsilon}\left(Z_{0}\right)}
	\frac{\omega_{1} \left(\left|S-Z_{\ast}\right|\right)}{\left|S-Z_{\ast}\right|}\, |d S| \le \]\
	\[\le \frac{c}{\left(\mathcal{R}\left(Z_{0}\right)\right)
		^{\beta} \omega_{1}\left(({3}/{8})\,\mathcal{R}\left(Z_{0}\right)\right)}
	\int\limits_{0}^{4\varepsilon} \frac{\omega_{1} \left(\eta \right)}{\eta}\, d\eta.\]
	
	Аналогічно отримуємо оцінку
	\[\left|i_{2}\right| \le \frac{c}{\left(\mathcal{R}\left(Z_{0}\right)\right)
		^{\beta} \omega_{1}\left(({1}/{2})\,\mathcal{R}\left(Z_{0}\right)\right)}
	\int\limits_{0}^{2\varepsilon} \frac{\omega_{1} \left(\eta \right)}{\eta}\, d\eta.\]

	З урахуванням співвідношень \eqref{omeg-locBound-circle}, \eqref{omeg-locBound-diff}
	і нерівностей $|Z - Z_{\ast}| \le  \varepsilon$,   $  |S-Z| \ge (1/2)\left|S-Z_0\right| $, $\mathcal{R}(S)\ge (1/2)\mathcal{R}(Z_0)$, які справджуються  для кожного
	$S \in  \Gamma \setminus \left(\Gamma_{\frac{1}{2} \mathcal{R}(Z_0)}\left\langle \mathcal{X}\right\rangle \, \bigcup \, \Gamma_{2 \varepsilon}\left(Z_{0}\right)\right)$, отримуємо оцінки інтегралів $i_{3}$ і $i_{4}$:
	\[ \left|i_{3}\right| \le c\, \frac{\left|Z_{\ast} -Z_{0}\right|}{\left(\mathcal{R}\left(Z_{0}\right)\right)^{\beta+\beta_0} \omega_{2}\left(({1}/{2})\,\mathcal{R}\left(Z_{0}\right)\right)}
	\int\limits_{ \Gamma \setminus \left(\Gamma_{\frac{1}{2} \mathcal{R}(Z_0)}\left\langle \mathcal{X}\right\rangle \, \bigcup \, \Gamma_{2 \varepsilon}\left(Z_{0}\right)\right)}\frac{\omega_2 \left(|S-Z_{0}|\right)}
	{\left|S-Z\right|\left|S-Z_{0}\right|}\, |d S| \le \]
	\[\le \frac{c}{\left( \mathcal{R}\left(Z_0\right)\right)^{\beta+\beta_0}\omega_{2} \left((1/2)\mathcal{R}\left(Z_0\right)\right)} \, \varepsilon \int\limits_{ 2\varepsilon}^{2}\frac{\omega_2 (\eta)}{\eta^2}\, d\eta ,\]
	\[\left|i_{4}\right| \le c \, |Z-Z_0|\int\limits_{ \Gamma \setminus \left(\Gamma_{\frac{1}{2} \mathcal{R}(Z_0)}\left\langle \mathcal{X}\right\rangle \, \bigcup \, \Gamma_{2 \varepsilon}\left(Z_{0}\right)\right)}\frac{\omega_1 \left(|S-Z_0|\right)}{\left( \mathcal{R}\left(S\right)\right)^{\beta} \omega_{1}\left((1/2)\mathcal{R}\left(Z_0\right)\right)  |S-Z_0| |S-Z|}|dS| \le \]
	\[ \le \frac{c}{\left( \mathcal{R}\left(Z_0\right)\right)^{\beta}\omega_{1}\left((1/2)\mathcal{R}\left(Z_0\right)\right)} \, \varepsilon \int\limits_{2 \varepsilon}^{2}\frac{\omega_1(\eta)}{\eta^2}\, d\eta. \]

	Використовуючи наведені оцінки і спрямовуючи $\varepsilon$ до нуля, отримуємо рівність (\ref{limChOmeg-CP}).
	Лему доведено.
	
	\medskip
	\begin{lemma}\label{BoDiffBikhIntCorPs}
		{\it Нехай  функція $\Omega \colon \left(\Gamma \setminus \mathcal{X}\right) \times \left(\Gamma \setminus \mathcal{X}\right) \longrightarrow \mathbb{C} $ задовольняє оцінку \eqref{omeg-locBound-circle}
			і оцінку
			\begin{multline} \label{omeg-locBound-diffBound}
				\left| \Omega\left(S,T\right)- \Omega\left(S,T_0 \right) \right|
				\le \frac{c  }{\left(\mathcal{R}(S)\right)^{\beta} \left(\mathcal{R}(T_0)\right)^{\beta_0}}
				\min\left\{\frac{\omega_{2}\left(|S-T_0|\right)}
				{\omega_{2}\left(\frac{1}{2}\mathcal{R}(T_0)\right)}, \, 1 \right\}
				\frac{|T-T_0|}{|S-T_0|}\\[2mm]
				\forall T_0 \in \Gamma \setminus \mathcal{X} \quad \forall T \in \Gamma_{\frac{1}{4}\mathcal{R}(T_0)}(T_0) \quad
				\forall S \in \Gamma \setminus \mathcal{X}: |S-T_0| \ge 2|T-T_0|,
			\end{multline}
			де стала $c$ не залежить від $S$, $T$, $T_0$\emph{;} $\beta \in (0; 1)\emph{;} \, \beta_0 \ge 0$\emph{;}
			$\omega_k \colon (0; +\infty) \longrightarrow (0; +\infty)$ при $k=1,2$,~---
			неспадні обмежені  функції такі, що
			\[\omega_2(\varepsilon)\to 0, \quad \varepsilon \to 0, \]
			і функція $\omega_1$ задовольняє умову Діні \eqref{dini}.
			
			Тоді для довільного $r<\frac{1}{4}\mathcal{R}_0$ і довільних $T, T_0 \in \Gamma\setminus  \Gamma_{r}\left\langle \mathcal{X}\right\rangle$ таких, що $|T-T_0|=\varepsilon \le \frac{1}{8}r$, справедлива оцінка
			\begin{multline}
				\label{DiffBoudIntCo}
				\left| \int\limits_{\Gamma}\frac{\Omega(S,T)}{S-T}\, dS-  \int\limits_{\Gamma}\frac{\Omega(S,Z_{0})}{S-Z_{0}}\, dS\right|
				\le \frac{c}{r^{\beta}} \left(\frac{\varepsilon}{r^{1+\beta_0}} + \frac{\varepsilon}{r^{\beta_0}\,\omega_{2}((1/2)r)} \int\limits_{\varepsilon}^{2} \frac{\omega_2\left(\eta\right)}{\eta^2}\, d\eta + \right.\\ \left.+\frac{1}{\omega_{1}((1/2)r)}\left(\int\limits_{0}^{\varepsilon}\frac{\omega_{1}\left(\eta\right)}{\eta}\, d\eta +\varepsilon
				\int\limits_{\varepsilon}^{2} \frac{\omega_1\left(\eta\right)}{\eta^2}\, d\eta \right)
				\right),
			\end{multline}
			де стала $c$ не залежить від $r$, $T$ і $T_0$.
		}
	\end{lemma}
	
	{\bf \em Доведення.}
	Нехай $T$, $T_0  \in \Gamma\setminus  \Gamma_{r}\left\langle \mathcal{X}\right\rangle$ такі, що $|T-T_0|=\varepsilon \le (1/8)r$. Тоді
	\[ \int\limits_{\Gamma}\frac{\Omega(S,Z)}{S-T}\, dS-  \int\limits_{\Gamma}\frac{\Omega(S,T_{0})}{S-T_{0}}\, dS=
	\int\limits_
	{\Gamma_{\frac{1}{2} r}\left\langle \mathcal{X}\right\rangle }
	\frac{\Omega\left(S,T\right)- \Omega\left(S,T_{0}\right)}{S-T}\,d S +\]\[+  \int\limits_{\Gamma_{\frac{1}{2} r}\left\langle \mathcal{X}\right\rangle}
	\Omega(S,T_{0})\left(\frac{1}{S-T}- \frac{1}{S- T_{0}}\right)\, dS+ \int\limits_{\Gamma_{2\varepsilon}(T_0)}\frac{\Omega(S,T)}{S-T}\, dS-\int\limits_{\Gamma_{2\varepsilon}(T_0)}\frac{\Omega(S,T_0)}{S-T_0}\, dS+\]
	
	\[+\int\limits_{\Gamma  \setminus\,
		\Gamma_{\frac{1}{2} r}\left\langle \mathcal{X}\right\rangle\setminus \Gamma_{2\varepsilon}(T_0)}
	\frac{\Omega\left(S,T\right)- \Omega\left(S,T_{0}\right)}{S-T}\,d S + \]
	\[+
	\int\limits_{\Gamma  \setminus\,
		\Gamma_{\frac{1}{2} r}\left\langle \mathcal{X}\right\rangle\setminus \Gamma_{2\varepsilon}(T_0)}
	\Omega(S,T_{0})\left(\frac{1}{S-T}- \frac{1}{S- T_{0}}\right)\, dS
	=:  J_{1} + J_{2} + J_{3} - J_{4}+J_{5}+J_6.\]

Далі у доведенні через\, $c$\, позначено сталі, значення яких не залежить  від $S$, $T$ і $T_0$,
	але, взагалі кажучи, різні навіть у межах одного ланцюжка нерівностей.

	Інтеграли $J_1$, $J_2$ оцінюються аналогічно до інтегралів $j_1$, $j_2$ в лемі~\ref{contBikhIntCorPs}:
	\[\left|J_1\right|+ \left|J_2\right| \le c \frac{\varepsilon}{r^{2+\beta_0}}  \int\limits_
	{\Gamma_{\frac{1}{2} r}\left\langle \mathcal{X}\right\rangle}  \frac{|dS|}{\left(\mathcal{R}(S)\right)^{\beta}}  \le c\, \frac{\varepsilon}{r^{1+\beta+\beta_0}}.\]
	
	Інтеграли $J_3$, $J_4$, $J_5$, $J_6$ оцінюються аналогічно відповідно до інтегралів $i_1$, $i_2$, $i_3$, $i_4$ в лемі~\ref{contBikhIntCorPs}:
	\[ \left|J_3\right|+ \left|J_4\right| \le c \, \frac{1}{r^{\beta} \, \omega_1\left(\mathcal{R}(T)/2\right) }\int\limits_{\Gamma_{2 \varepsilon}(T_0)} \frac{\omega_{1}\left(|S-T|\right)}{|S-T|} \,|dS| +\]
	\[+
	\frac{1}{r^{\beta} \, \omega_1\left(\mathcal{R}(T_0)/2\right) }\int\limits_{\Gamma_{2 \varepsilon}(T_0)} \frac{\omega_{1}\left(|S-T_0|\right)}{|S-T_0|} \,|dS|\le  \frac{c}{r^{\beta}\, \omega_1 \left((1/2)r\right)}
	\int\limits_{0}^{3\varepsilon} \frac{\omega_{1}\left(\eta\right)}{\eta} \,d \eta .\]
	
	\[\left|J_5\right| \le
	\frac{c}{r^{\beta+\beta_0} \omega_2 \left((1/2)\mathcal{R}
		\left(T_0\right)\right)}\,\,
	\varepsilon
	\int\limits_{ \Gamma \setminus \Gamma_{\frac{1}{2} r}\left\langle \mathcal{X}\right\rangle \setminus \, \Gamma_{2 \varepsilon}\left(T_{0}\right)}
	\frac{\omega_2 \left(|S-T_0|\right)}
	{  |S-T_0| |S-T|} |dS| \le  \]
	\[\le \frac{c}{r^{\beta+\beta_0}\omega_2 \left((1/2)r\right)}\, \varepsilon \int\limits_{2\varepsilon}^{2} \frac{\omega_{2}\left(\eta\right)}{\eta^2}\, d\eta.\]
	
	\[\left|J_6\right| \le
	\frac{c}{r^{\beta} \omega_1 \left((1/2)\mathcal{R}
		\left(T_0\right)\right)}\,\,
	\varepsilon
	\int\limits_{ \Gamma \setminus \Gamma_{\frac{1}{2} r}\left\langle \mathcal{X}\right\rangle \setminus \, \Gamma_{2 \varepsilon}\left(T_{0}\right)}
	\frac{\omega_1 \left(|S-T_0|\right)}
	{  |S-T_0| |S-T|} |dS| \le  \]
	\[\le \frac{c}{r^{\beta}\omega_1 \left((1/2) r\right)}\, \varepsilon \int\limits_{2\varepsilon}^{2} \frac{\omega_{1}\left(\eta\right)}{\eta^2}\, d\eta.\]
	
	Лему доведено.

	\medskip
	\begin{lemma}
		\label{IneqOmegMoDiffPract}
		{\it Нехай існує $r_0<\frac{1}{4}\mathcal{R}_0$ таке, що
			функція $\Omega \colon \left(\Gamma \setminus \mathcal{X}\right) \times \left(\overline{\mathcal{U}} \setminus \mathcal{X}\right) \longrightarrow \mathbb{C} $
			при кожному  $X_j$, $j=\overline{1,m}$,
			задовольняє оцінки
			\begin{equation}
				\label{ineqOmonDomCp}
				\left|\Omega\left(S,Z\right) \right| \le
				c\, \frac{\left(\max\left\{
					\mathcal{R}(S),\mathcal{R}(Z)  \right\}\right)^{\gamma_j'}}{(\mathcal{R}(S))^{\gamma_j + \gamma_j'}}  \qquad  \forall S \in \Gamma \setminus \mathcal{X} \,\,\,
				\forall Z \in \mathcal{U}_{r_0}\left(X_j\right), 
			\end{equation}
			\begin{multline}
				\label{ineqOmonBoundCp}
				\left|\Omega\left(S,T\right) \right| 
				\le c\, \frac{\left(\max\left\{
					\mathcal{R}(S),\mathcal{R}(T)  \right\}\right)^{\gamma_j'}}{(\mathcal{R}(S))^{\gamma_j + \gamma_j'}}
				\min\left\{\frac{\omega_1^{\ast}\left(|S-T|\right)}{\omega_1^{\ast}\left(\frac{1}{2} \mathcal{R}(T)\right)}, 1\right\} \qquad  \forall S \in \Gamma \setminus \mathcal{X} \\ 
				\forall T \in \Gamma_{r_0}\left(X_j\right)\setminus\{X_j\}, 
			\end{multline}
			\begin{multline}
				\label{ineqDiffOmBound}
				\shoveleft{
					\left|\Omega(S,Z)-\Omega(S,Z_0)\right| \le  c\, \frac{1}{\left(\mathcal{R}(Z_0)\right)^{\gamma_j}}\frac{|Z-Z_0|}{|S-Z_0|}} \qquad
				\forall Z_0 \in
				\Gamma_{r_0}\left(X_j\right)\setminus\{X_j\}\,\,\, 
				\forall Z \in \overline{\mathcal{U}}_{\frac{1}{4}\mathcal{R}(Z_0)}(Z_0)\\[1mm] 
				\forall S \in \Gamma_{\frac{1}{2}\mathcal{R}(Z_0)}(Z_0) : |S-Z_0| \ge 2|Z-Z_0|, 
			\end{multline} 
			де стала
			$c$ не залежить від    $S$, $Z$, $T$ і $Z_0$; додатні числа  $\gamma_j>0$ та $\gamma_j'>0$ задовольняють умову  $\gamma_j+ \gamma_j'<1$; 
			$\omega_1^{\ast} \colon (0; +\infty) \longrightarrow (0; +\infty)$~--- неспадна обмежена  функція, яка задовольняє умову
			\begin{equation}\label{dini-om1Pr}
				\int\limits_{0}^{\varepsilon}\frac{\omega_1^{\ast}(\eta)}{\eta}\, d\eta < c\, \omega_1^{\ast}(\varepsilon)  \qquad \forall \varepsilon \in  \left(0; \frac{1}{2} r_0 \right),
			\end{equation}
			де $c$~--- деяка стала, яка не задежить від $\varepsilon$;
			
			Тоді в околі кожної точки $X_j \in \mathcal{X}$, $j=\overline{1,m}$, виконується оцінка
			\begin{equation}
				\label{EstIntCtNearCps}
				\left|\int\limits_{\Gamma}
				\frac{\Omega(S,Z)}{S-Z}\, dS \right| \le \frac{c}{\left(\mathcal{R}(Z)\right)^{\gamma_j}} \qquad  \forall Z \in 	\mathcal{U}_{r_0/2}\left(X_j\right), 
			\end{equation}
			де стала $c$ не залежить від  $Z$. 
		}
	\end{lemma}
	
	{\bf \em Доведення.}
	Нехай $Z\in \mathcal{U}_{r_0/2}\left(X_j\right)$. Для доведення оцінки \eqref{EstIntCtNearCps} використаємо рівність
		\begin{multline*}
		\int\limits_{\Gamma}\frac{\Omega(S,Z)}{S-Z} dS=\int\limits_{\Gamma_{(7/8)\mathcal{R}(Z)}\left(X_j\right)} \frac{\Omega(S,Z)}{S-Z}\, dS +  \int\limits_{\Gamma_{(1/2)\mathcal{R}_0}\left(X_j\right)\setminus \Gamma_{(9/8)\mathcal{R}(Z)}\left(X_j\right)}
		\frac{\Omega(S,Z)}{S-Z} \, dS+\\[2mm]
		+\int\limits_{\Gamma_{(9/8)\mathcal{R}(Z)}\left(X_j\right) \setminus \, \Gamma_{(7/8)\mathcal{R}(Z)} \left(X_j\right)}\frac{\Omega(S,Z)}{S-Z} \, dS +  \int\limits_{\Gamma \setminus \, \Gamma_{(1/2)\mathcal{R}_0}\left(X_j\right)}
		\frac{\Omega(S,Z)}{S-Z} \, dS =: \sum _{k=1}^{4} I_{k}
	\end{multline*}
	і оцінимо інтеграли $I_{k}$,\, $k=\overline{1,4}$.

Скрізь у доведенні через\, $c$\, позначатимуться сталі, значення яких не залежить  від $S$, $Z$ і точки $Z_*$, яку буде введено далі,
	але, взагалі кажучи, різні навіть у межах одного ланцюжка нерівностей.

	З урахуванням нерівності \eqref{ineqOmonDomCp} та співвідношень $|S-Z| \ge (1/8)\mathcal{R}(Z)$, $\mathcal{R}(S) \le  (7/8) \mathcal{R}(Z)$, що виконуються для всіх $S \in \Gamma_{(7/8)\mathcal{R}(Z)}(X_j)$, маємо
	\[\left|I_{1}\right| \le \frac{8}{\mathcal{R}(Z)}\int\limits_{\Gamma_{(7/8)\mathcal{R}(Z)}\left(X_j\right)} \left|\Omega(S,Z)\right|\, \left| dS \right| \le
	\frac{c}{\left(\mathcal{R}(Z)\right)^{1-\gamma_j'}}
	\int\limits_{\Gamma_{(7/8)\mathcal{R}(Z)}\left(X_j\right)} \frac{|dS|}{\left(\mathcal{R}(S)\right)^{\gamma_j+\gamma_j'}} \le
	\]
	\[\le c \frac{1}{\left(\mathcal{R}(Z)\right)^{1-\gamma_j'}}\int\limits_{0}^{(7/8)\mathcal{R}(Z)}\frac{d\eta}{\eta^{\gamma_j+\gamma_j'}} \le \frac{c}{\left(\mathcal{R}(Z)\right)^{\gamma_j}}. \]
	
	Враховуючи співвідношення $|S-Z|\ge \left(1/9\right)\mathcal{R}(S)$ і $\mathcal{R}(S) \ge \mathcal{R}(Z)$, що виконуються для всіх $S \in \Gamma_{(1/2)\mathcal{R}_0}\left(X_j\right)\setminus \Gamma_{(9/8)\mathcal{R}(Z)}\left(X_j\right)$,
	та нерівність \eqref{ineqOmonDomCp},  отримуємо оцінку інтеграла $I_{2}$:
	\[\left|I_{2}\right| \le 9\,
	\int\limits_{\Gamma_{(1/2)\mathcal{R}_0}\left(X_j\right)\setminus \Gamma_{(9/8)\mathcal{R}(Z)}\left(X_j\right)}
	\frac{|\Omega(S,Z)|}{\mathcal{R}(S)}\,|dS|
	\le c \int\limits_{\Gamma_{(1/2)\mathcal{R}_0}\left(X_j\right)\setminus \Gamma_{(9/8)\mathcal{R}(Z)}\left(X_j\right)} \frac{|dS|}{\left(\mathcal{R}(S)\right)^{1+\gamma_j}}\le\]
	\[\le c \int\limits_{(9/8)\mathcal{R}(Z)}^{(1/2)\mathcal{R}_0}\frac{d\eta}{\eta^{1+\gamma_j}} \le  \frac{c}{\left(\mathcal{R}(Z)\right)^{\gamma_j}}.\]

	Враховуючи, нерівність \eqref{ineqOmonDomCp} і нерівність $|S-Z|>(1/4) \mathcal{R}_0$, що виконується для
	всіх  $S \in \Gamma \setminus \Gamma_{\frac{1}{2}\mathcal{R}_0}\left(X_j\right)$, отримуємо оцінку інтеграла $I_{4}$:
	\[\left|I_{4}\right| \le \frac{4}{\mathcal{R}_0}\int\limits_{\Gamma \setminus \Gamma_{\frac{1}{2}\mathcal{R}_0}\left(X_j\right)} \left|\Omega(S,Z)\right|\, \left| dS \right| \le\,
	c\,\sum\limits_{k=1}^{m}\int\limits_{\Gamma}\frac{|dS|}{\left(\mathcal{R}(S)\right)^{\gamma_k+\gamma_k'}} 
	\le\, c\,\sum\limits_{k=1}^{m}\int\limits_{0}^{2}\frac{d\eta}{\eta^{\gamma_k+\gamma_k'}} \le\, c\, \le \frac{c}{\left(\mathcal{R}(Z)\right)^{\gamma_j}}. \]
	
	Нехай $Z_{\ast}$ --- точка кола $\Gamma$, найближча до точки $Z$.
	Позначимо\, $\delta:= \left|Z-Z_{\ast}\right|\equiv 1- |Z|$\, і\, $\widetilde{\Gamma}:=\Gamma_{(9/8)\mathcal{R}(Z)}\left(X_j\right) \setminus \, \Gamma_{(7/8)\mathcal{R}(Z)}\left(X_j\right)$.
	
	Для оцінки інтеграла $I_3$ розглянемо два випадки.
	
	У випадку
	$\delta \ge (1/16) \mathcal{R}(Z)$,
	враховуючи нерівність $\left|S-Z\right| \ge (1/16)\mathcal{R}(Z)$,
	що виконується для всіх $S \in \widetilde{\Gamma}$,
	та нерівність \eqref{ineqOmonDomCp}, а також співвідношення
	\begin{equation}\label{podv-estim-1}
		(7/8)\mathcal{R}(Z) \le \mathcal{R}(S) \le (9/8)\mathcal{R}(Z)
		\qquad \forall   S\in\widetilde{\Gamma},
	\end{equation}
	отримуємо
	\[\left|I_{3}\right| \le \frac{16}{\mathcal{R}(Z)}\int\limits_{\widetilde{\Gamma}} \left|\Omega(S,Z)\right|\,\left| dS \right|
	\le \frac{c}{\left(\mathcal{R}(Z)\right)^{1+\gamma_j}} \int\limits_{\Gamma_{(9/8)\mathcal{R}(Z)}(X_j)}|dS| \le \frac{c}{\left(\mathcal{R}(Z)\right)^{\gamma_j}}.\]
	
	Розглянемо, нарешті, випадок $\delta < (1/16) \mathcal{R}(Z)$.
	Подамо $I_{3}$ у вигляді суми інтегралів:
		\[I_{3} = \int\limits_{\widetilde{\Gamma}} \frac{\Omega\left(S,Z\right)-\Omega\left(S,Z_{\ast}\right)}{S-Z}\, dS +
		\int\limits_{\widetilde{\Gamma}} \frac{\Omega\left(S,Z_{\ast}\right)}{S-Z}\, dS=: I_{3,1} + I_{3,2}.\]
		
	Інтеграл $I_{3,1}$ у свою чергу оцінимо сумою двох інтегралів:
	\[\left|I_{3,1}\right|\le  \int\limits_{\Gamma_{2\delta}\left(Z_{\ast}\right)}\frac{\left|\Omega\left(S,Z\right)-\Omega\left(S,Z_{\ast}\right) \right| }{|S - Z|}\, |dS|
	+\int\limits_{\widetilde{\Gamma} \setminus\Gamma_{2\delta}\left(Z_{\ast}\right)}\frac{\left|\Omega\left(S,Z\right)-\Omega\left(S,Z_{\ast}\right) \right|}{|S - Z | }\, |dS|=:
	I_{3,1}' + I_{3,1}''. \]
	
	Враховуючи нерівність $|S-Z|\ge  \delta$, що виконується для всіх $S \in \Gamma$,
	та нерівності \eqref{ineqOmonDomCp} і \eqref{ineqOmonBoundCp}, а також
	cпіввідношення \eqref{podv-estim-1} і
	\begin{equation}\label{podv-estim-2}
		(15/16)\mathcal{R}(Z) \le \mathcal{R}(Z_{\ast}) \le (17/16)\mathcal{R}(Z),
	\end{equation}
	отримуємо
	\[\left|I_{3,1}'\right| \le \int\limits_{\Gamma_{2\delta}\left(Z_{\ast}\right)}\frac{\left|\Omega\left(S,Z\right)\right| + \left|\Omega\left(S,Z_{\ast}\right) \right|}{\left| S - Z \right|}\, |dS| %\le\]
	%\[
	\le \frac{c}{\delta\,\left(\mathcal{R}(Z)\right)^{\gamma_j}}
	\int\limits_{\Gamma_{2\delta}\left(Z_{\ast}\right)}|dS| \le
	\frac{c}{\left(\mathcal{R}(Z)\right)^{\gamma_j}}. \]
	
	Враховуючи нерівність
	\begin{equation}\label{estim-znyzu-1}
		|S-Z|\ge (1/2)|S-Z_{\ast}|	\qquad \forall   S\in\Gamma
	\end{equation}
	і нерівність \eqref{ineqDiffOmBound}, а також cпіввідношення \eqref{podv-estim-2}, отримуємо
	\[ \left|I_{3,1}''\right| \le 2\,
	\int\limits_{\widetilde{\Gamma} \setminus\Gamma_{2\delta}\left(Z_{\ast}\right)}\frac{\left|\Omega\left(S,Z\right)-\Omega\left(S,Z_{\ast}\right) \right|}{|S - Z_{\ast} | }\, |dS| \le \]
	\[\le
	c\, \frac{\delta}{\left(\mathcal{R}\left(Z_{\ast}\right)\right)^{\gamma_j}}
	\int\limits_{\widetilde{\Gamma} \setminus\Gamma_{2\delta}} \frac{|dS|}{\left|S-Z_{\ast}\right|^2} \le
	c\, \frac{\delta}{\left(\mathcal{R}\left(Z \right)\right)^{\gamma_j}}
	\int\limits_{2\delta}^{2}\frac{d\eta}{\eta^2} \le
	\frac{c}{\left(\mathcal{R}\left(Z \right)\right)^{\gamma_j}}.
	\]

	При оцінці інтеграла $I_{3,2}$ спочатку використовуємо нерівність \eqref{estim-znyzu-1}, після чого
	оцінюємо його сумою двох інтегралів:
	\[\left|I_{3,2} \right| \le
	\int\limits_{\widetilde{\Gamma}\cap\Gamma_{(1/2)\mathcal{R}(Z_{\ast})}(Z_{\ast})} \frac{2\,\left|\Omega\left(S,Z_{\ast}\right)\right|}{\left|S- Z_{\ast} \right|} \, |dS|
	+ \int\limits_{\widetilde{\Gamma}\setminus\Gamma_{(1/2)\mathcal{R}(Z_{\ast})}(Z_{\ast})} \frac{2\,\left|\Omega\left(S,Z_{\ast}\right)\right|}{\left|S- Z_{\ast} \right|} \, |dS|=:
	I_{3,2}' + I_{3,2}''. \]
	
	Враховуючи нерівність  \eqref{ineqOmonBoundCp},
	cпіввідношення \eqref{podv-estim-1} і \eqref{podv-estim-2},
	а також умову \eqref{dini-om1Pr}, отримуємо:
	\[|I_{3,2}'|\le \frac{c}{\left(\mathcal{R}(Z)\right)^{\gamma_j}\omega_{1}^{\ast}\left((1/2)\mathcal{R}(Z_{\ast})\right)} \int\limits_{\Gamma_{(1/2)\mathcal{R}(Z_{\ast})}(Z_{\ast})}\frac{\omega_{1}^{\ast}\left( \left|S-Z_{\ast} \right|\right)
	}{\left|S-Z_{\ast} \right|}\, |dS| \le \]
	\[ \le \frac{c}{\left(\mathcal{R}(Z)\right)^{\gamma_j}\omega_{1}^{\ast}\left((1/2)\mathcal{R}(Z_{\ast})\right)}
	\int\limits_{0}^{\frac{1}{2}  \mathcal{R}\left(Z_{\ast}\right)} \frac{\omega_{1}^{\ast}\left( \eta \right)
	}{\eta}\, d \eta \le
	\frac{c}{\left(\mathcal{R}(Z)\right)^{\gamma_j}}.\]
	
	Нарешті, враховуючи нерівність  \eqref{ineqOmonDomCp} та cпіввідношення \eqref{podv-estim-1} і \eqref{podv-estim-2},
	отримуємо:
	\[|I_{3,2}''|\le \frac{4}{\mathcal{R}(Z_{\ast})}\int\limits_{\widetilde{\Gamma}\setminus\Gamma_{(1/2)\mathcal{R}(Z_{\ast})}(Z_{\ast})} \left|\Omega(S,Z_{\ast})\right|\,\left| dS \right|
	\le \frac{c}{\left(\mathcal{R}(Z)\right)^{1+\gamma_j}} \int\limits_{\Gamma_{(9/8)\mathcal{R}(Z)}(X_j)}|dS| \le \frac{c}{\left(\mathcal{R}(Z)\right)^{\gamma_j}}.\]
	
	Лему доведено.

%	\subsection{Деякі спеціальні властивості конформного відображення круга на область з кусково-гладкою межею}
{\bf 3.2. Деякі спеціальні властивості конформного відображення круга на область з кусково-гладкою межею.}
		Скрізь надалі припускаємо існування неперервної  контурної похідної
	\[\sigma'(T):= \lim\limits_{S  \to T,\,
		S\in \Gamma \setminus \mathcal{X}}\frac{\sigma(S)-\sigma(T)}{S-T} \qquad \forall\, T \in \Gamma \setminus \mathcal{X},\]  конформного відображення
	$\sigma(Z)$, яка, крім того,  відмінна від нуля  на $\Gamma \setminus \mathcal{X}$. %:

	З результатів роботи \cite{Walsh-Sewell-1940} або більш загальної
	теореми~2.7.6 \cite{tamrazov} випливає, що в усіх точках множини $\Gamma \setminus \mathcal{X}$ існує наступна тілесна похідна, значення якої співпадають зі значеннями контурної похідної, тобто виконується рівність
	\[  \lim\limits_{Z  \to Z_0,\,
		Z\in \overline{\mathcal{U}}}\frac{\sigma(Z)-\sigma(Z_0)}{Z-Z_0}=\sigma'\left(Z_0\right) \qquad \forall\, Z_0 \in \Gamma \setminus \mathcal{X},\]
	і, крім того, $\sigma'\left(Z\right)$ є неперервною при всіх  $Z \in \overline{\mathcal{U}} \setminus \mathcal{X}$.
	
	Введемо позначення 
	\[\sigma_{1}\left(Z\right):= \mathrm{Re}\, \sigma\left(Z\right),\qquad
	\sigma_{2}\left(Z\right):= \mathrm{Im}\, \sigma\left(Z\right), \]
	\[ d\left(S,Z\right):=
	\frac{\sigma\left(S\right) -\sigma\left(Z\right)}{S-Z},
	\qquad
	d_{k}\left(S,Z\right):=\frac{\sigma_{k}\left(S\right) -\sigma_{k}\left(Z\right)}{S-Z},\,\,\, k=1,2, \]
	\[\forall\,Z \in \overline{\mathcal{U}} \setminus \mathcal{X}\quad \forall\,S \in \Gamma \setminus \mathcal{X}: Z\ne S. \]
	
	Для  $r_0 \in (0;\mathcal{R}_0/4)$ введемо в розгляд криву
	\[\gamma_{r_0}:=\left( \Gamma \setminus
	\Gamma_{r_0}\langle{\mathcal{X}\rangle} \right)
	\bigcup
	\left(\,\bigcup_{j=1}^{m}\{S\in \overline{\mathcal{U}}: \left| S-X_{j}\right| = r_0\}\right),
	\]
	орієнтація якої індукована орієнтацією кола $\Gamma$.
	Позначимо через $\mathcal{U}_{r_0}$  область, обмежену кривою $\gamma_{r_0}$, а через $\overline{\mathcal{U}_{r_0}}$ --- замикання цієї області.
	
	Оскільки $\sigma'\left(T\right)\ne 0$   для всіх $T \in \overline{\mathcal{U}_{r_0}}$,
	то	існують додатні сталі   $m_{\ast}(r_0)$ і $M_{\ast}(r_0)$, залежні від $r_0$, такі, що виконується подвійна нерівність
	\begin{equation}
		\label{sigPrle0Xr}
		m_{\ast}(r_0)  \le \left| \sigma'\left(T\right)\right| \le M_{\ast}(r_0) \qquad  \forall\, T \in \gamma_{r_0}.
	\end{equation}

	Розглянемо параметричне представлення $S=S(s)$ кривої $\gamma_{r_0}$, де $s \in [0,l(\gamma_{r_0})]$ --- натуральний параметр і $l(\gamma_{r_0})$~--- довжина кривої  $\gamma_{r_0}$, при цьому $S(0)=S(l(\gamma_{r_0}))$.
	Для кожної точки  $T\in\gamma_{r_0}$, відмінної від точки $S(0)$, відповідною малою літерою $t$ позначатимемо відповідне значення натурального параметра, при якому виконується рівність  $T=S(t)$.
	
	Для пари різних точок $S,T\in \gamma_{r_0}$ через $\gamma_{r_0}[S,T]$ позначимо ту дугу кривої $\gamma_{r_0}$ з кінцями $S$ і $T$,
	яка має небільшу довжину $l\left(\gamma_{r_0}[S,T]\right)$. Очевидно, що у випадку, коли $S(0)\notin\gamma_{r_0}[S,T]$, виконується рівність
	$l\left(\gamma_{r_0}[S,T]\right)=|s-t|$.

	Очевидно, що існує абсолютна стала $K>0$ така, що виконується подвійна нерівність
	\begin{equation}
		\label{duga-khorda}
		1 \le \frac{l\left(\gamma_{r_0}[S,T]\right)}{|S-T|}\le K \qquad \forall\,S,T\in \gamma_{r_0} : S\ne T.
	\end{equation}

	\medskip
	\begin{lemma}
		\label{DbiggerZe}
		{\it Нехай  контурна похідна $\sigma'$  неперервна та не обертається в   нуль на  $\Gamma \setminus \mathcal{X}$.
			Тоді для кожного  $r_0 \in (0;\mathcal{R}_0/4)$ існує додатна стала 
			$M(r_0)$,   залежна від $r_0$, така,
			що при $k=1,2$ виконується нерівність}
		\begin{equation}
			\label{dTypPokhInUbezX}
		\left|d_k\left(S,Z\right) \right| \le M(r_0) \qquad \forall S\in \gamma_{r_0} \quad \forall Z \in \overline{\mathcal{U}_{r_0}} : S \ne Z.
		\end{equation}
	\end{lemma}

	{\bf \em Доведення.}
	При доведенні оцінки \eqref{dTypPokhInUbezX} розглянемо спочатку випадок $Z=T$ і  $S, T \in \gamma_{r_0} : S \ne T$.
	Не зменшуючи загальності міркувань, можна вважати, що  $S(0)\notin\gamma_{r_0}[S,T]$.

	З урахуванням верхніх оцінок у подвійних нерівностях   \eqref{sigPrle0Xr} і \eqref{duga-khorda} при $k=1,2$ маємо співвідношення
	\begin{multline*}
		|d_k(S,T)|\le |d(S,T)|= \frac{1}{|S-T|}\biggl|\,\int\limits_{\gamma_{r_0}[S,T]} \sigma'(S)\, dS\biggr|
		\le \frac{1}{|S-T|}\int\limits_{\gamma_{r_0}[S,T]} \bigl|\sigma'(S)\bigr|\,|dS| \le \\
		\le M_{\ast}(r_0)\,\frac{l\left(\gamma_{r_0}[S,T]\right)}{|S-T|}
		\le K M_{\ast}(r_0).
	\end{multline*}
	
	Отже, доведено верхню оцінку
	\begin{equation} \label{Ver-oz-cont}
		|d(S,T)|\le M(r_0) \qquad \forall\,S,T\in\gamma_{r_0} : S\ne T
	\end{equation}
	при     $M(r_0):= K M_{\ast}(r_0)$.

	Залишається довести, що нерівність \eqref{Ver-oz-cont} залишається в силі при заміні
	змінної $T\in\gamma_{r_0}$ на змінну $Z\in \overline{\mathcal{U}_{r_0}}$.
	Для цього зауважимо, що для кожної фіксованої точки $S\in  \gamma_{r_0}$ функція $d(S,Z)$ за змінною $Z$
	є голоморфною в облаcті  $\mathcal{U}_{r_0}$ і, крім того, неперервено продовжується  на  $\gamma_{r_0}\setminus\{S\}$.
	Неперервна продовжуваність цієї функції у точку $Z=S$ випливає з існування тілесної похідної відображення $\sigma$ у цій точці.
	Отже, для кожної фіксованої точки $S\in  \gamma_{r_0}$ функція $d(S,Z)$ за змінною $Z$ є неперервною в $\overline{\mathcal{U}_{r_0}}$ і голоморфною в області $ \mathcal{U}_{r_0}$.
	
	Тепер, спираючись на принцип максимума для голоморфних функцій, робимо висновок  про те, що нерівність \eqref{Ver-oz-cont} залишається в силі при заміні змінної $T\in\gamma_{r_0}$ на змінну $Z\in \overline{\mathcal{U}_{r_0}}$.
	
	Лему доведено.

	\medskip
	\begin{lemma}
		\label{Dbig-oz-znyzu}
		{\it Нехай  контурна похідна $\sigma'$  неперервна та не обертається в   нуль на  $\Gamma \setminus \mathcal{X}$.
			Тоді для кожного  $r_0 \in (0;\mathcal{R}_0/4)$ існує додатна стала $m(r_0)$,   залежна від $r_0$, така,
			що виконується нерівність}
		\begin{equation}
			\label{d-oz-znyzu}
			\left|d\left(S,Z\right) \right| \ge m(r_0) \qquad \forall S\in \gamma_{r_0} \quad \forall Z \in \overline{\mathcal{U}_{r_0}} : S \ne Z.
		\end{equation}
	\end{lemma}

	{\bf \em Доведення.}
	При доведенні оцінки \eqref{d-oz-znyzu} розглянемо спочатку випадок $Z=T$ і  $S, T \in \gamma_{r_0} : S \ne T$.
	Не зменшуючи загальності міркувань, можна вважати, що $S(0)\notin\gamma_{r_0}[S,T]$ і в точці $S(0)$ існує дотична до кривої $\gamma_{r_0}$.

	З існування контурної похідної $\sigma'(T)$ в усіх точках $T\in \gamma_{r_0}$ випливає існування односторонніх границь
	\[\sigma_{k}'(T_0+0):=\lim\limits_{t\to t_0+0} \sigma_{k}'(S(t)), \qquad k=1,2, \]
	і
	\[\sigma_{k}'(T_0-0):=\lim\limits_{t\to t_0-0} \sigma_{k}'(S(t)), \qquad k=1,2, \]
	для всіх $T_0=S(t_0)$, $t_0\in (0,l(\gamma_{r_0}))$,
	при цьому в кожній точці $T\in \gamma_{r_0}\setminus\{S(0)\}$, в якій існує дотична до кривої $\gamma_{r_0}$, виконуються рівності
	\[\sigma_{k}'(T+0)=\sigma_{k}'(T-0), \qquad k=1,2. \]
	
	З \eqref{sigPrle0Xr} випливає, що
	для кожної точки $T\in \gamma_{r_0}$, в якій існує дотична до кривої $\gamma_{r_0}$, виконується хоча б одна з двох нерівностей
	\begin{equation}\label{sigpNe0}
		\left|\sigma_{1}'\left(T\right) \right| \ge  m_{\ast}(r_0)/2 \quad \mbox{або} \quad
		\left|\sigma_{2}'\left(T\right) \right| \ge m_{\ast}(r_0)/2,
	\end{equation}
	а в кожній кутовій точці $T_0$ 
	кривої $\gamma_{r_0}$ виконується хоча б одна з двох пар нерівностей
	\begin{equation}\label{sigpNe0-kut}
		\left|\sigma_{1}'(T_0\pm 0) \right| \ge  m_{\ast}(r_0)/2 \quad \mbox{або} \quad
		\left|\sigma_{2}'(T_0\pm 0) \right| \ge m_{\ast}(r_0)/2.
	\end{equation}

	Оскільки функція $\sigma'$ неперервна на $\gamma_{r_0}$, то
	з нерівностей \eqref{sigpNe0}, \eqref{sigpNe0-kut} випливає, що для будь-якого $T \in \gamma_{r_0}$ існує $\delta(T)$
	таке, що $0<\delta(T)<r_0/2$ і для усіх $S\in \gamma_{r_0}$, які задовольняють умову $|S-T|<\delta(T)$, виконується хоча б одна з двох пар нерівностей
	\begin{equation}
		\label{siPravnOts}
		\left|\sigma_{1}'(S\pm 0)\right|> m_{\ast}(r_0)/3 \quad \mbox{aбо} \quad
		\left|\sigma_{2}'(S\pm 0)\right|> m_{\ast}(r_0)/3.
	\end{equation}
	
	Для кожного $k\in \{1,2\}$ і $t \in (0,l(\gamma_{r_0}))$
	розглянемо функцію
	\[\widetilde\sigma_{k}(t):=\sigma_{k}\left(T\right),\qquad \mbox{де} \quad T=S(t),\]
	для якої майже скрізь на $\gamma_{r_0}$ виконується очевидна рівність
	\begin{equation}
		\label{rivn-dlia-pokh}
		\widetilde\sigma_{k}'(t)=  \sigma_k'(T)S'(t).
	\end{equation}
	Введемо також у розгляд середні значення похідних $\widetilde\sigma_{k}'$, $k\in \{1,2\}$, на
	відрізку, що є прообразом дуги $\gamma_{r_0}[S,T]$ при відображенні $T=S(t)$:
	\begin{equation}
		\label{mean}
		\widetilde\sigma_{k,{\rm mean}}':=\frac{1}{s-t}\,\int\limits_{t}^{s} \widetilde\sigma_k'(s)\, ds.
	\end{equation}

	Оскільки $\gamma_{r_0}$~--- компакт, то існує скінченне покриття кривої $\gamma_{r_0}$ відкритими кругами
	$K_{\delta(T_j)}(T_j):=\{Z\in \mathbb{C}: |Z-T_j|< \delta(T_j)\}$, де $T_j\in\gamma_{r_0}$, $j=\overline{1,n}$.
	Для кожного  $j\in\{1,2, \dots, n\}$ позначимо через $Z_{j}'$ і
	$Z_{j}''$~--- точки перетину кола
	$\partial K_{\delta(T_j)}
	\left(T_j\right):=
	\left\{ Z\in \mathbb{C}: \left| Z- T_j \right|= \delta_j(T_j)
	\right\}$ та кривої  $\gamma_{r_0}$.
	Позначимо
	\[ \delta_0:= \min\limits_{j\ne l}\left\{\left| Z_{j}'-Z_{l}'\right|, \left| Z_{j}'-Z_{l}''\right|,  \left| Z_{j}''-Z_{l}''\right| \right\}.\]

	Якщо $|S-T|<\delta_0$, то існує  $l\in \{1,2, \dots, n\}$ таке, що
	$\gamma_{r_0} [S,T] \subset K_{\delta_l}\left(T_l\right)$. Тоді з урахуванням рівностей \eqref{rivn-dlia-pokh} і \eqref{mean},
	а також нижньої оцінки у подвійній нерівності \eqref{duga-khorda}, отримуємо співвідношення
		\begin{multline*}
		|d(S,T)|=\frac{1}{|S-T|} \biggl|\,\int\limits_{\gamma_{r_0}[S,T]} \left(\sigma_1'(S)+i\sigma_2'(S)\right)\, dS\biggr|=
		\frac{1}{|S-T|} \biggl|\,\int\limits_{t}^{s} \left(\widetilde\sigma_1'(s)+i\widetilde\sigma_2'(s)\right)\, ds\biggr|=\\
		=\Bigl|\widetilde\sigma_{1,{\rm mean}}' +
		i\widetilde\sigma_{2,{\rm mean}}' \Bigr|\,\frac{l\left(\gamma_{r_0}[S,T]\right)}{|S-T|}\ge \max_{k\in\{1,2\}}\min_{H\in\gamma_{r_0}[S,T]}\,\left|\sigma_{k}'(H\pm 0)\right|.
	\end{multline*}
	
	Оскільки $\gamma_{r_0} [S,T] \subset K_{\delta_l}\left(T_l\right)$, то наслідком співвідношень \eqref{siPravnOts} є нерівність
	\[ \max_{k\in\{1,2\}}\min_{H\in\gamma_{r_0}[S,T]}\,\left|\sigma_{k}'(H\pm 0)\right|\ge m_{\ast}(r_0)/3. \]
	
	Отже, за умови $|S-T|<\delta_0$ встановлено оцінку
	\[	|d(S,T)|\ge  m_{\ast}(r_0)/3. \]

	Якщо  $|S-T|\ge \delta_0$, то, враховуючи, що відображення  $\sigma$ є гомеоморфізмом на $\gamma_{r_0}$, можна стверджувати, що існує додатна стала $m_0(\delta_0)$, залежна від $\delta_0$, така, що виконується нерівність  $\Bigl|\sigma(S)-\sigma(T)\Bigr|\ge m_{0}(\delta_0)$ для всіх $S,T\in\gamma_{r_0} : |S-T|\ge \delta_0$. Тому при $|S-T|\ge \delta_0$ отримуємо оцінку
	\begin{equation} \label{Nyzh-oz-homeomorph}
		|d(S,T)|= \frac{\Bigl|\sigma(S)-\sigma(T)\Bigr|}{|S-T|} \ge m_{0}(\delta_0)/2.
	\end{equation}
	
	Отже, доведено нижню оцінку
	\begin{equation} \label{Nyzh-oz-cont}
		|d(S,T)|\ge m(r_0) \qquad \forall\,S,T\in\gamma_{r_0} : S\ne T
	\end{equation}
	при     $m(r_0):= \min\left\{m_{\ast}(r_0)/3, m_0(\delta_0)/2\right\}$.

	Залишається довести, що нерівність \eqref{Nyzh-oz-cont} залишається в силі при заміні
	змінної $T\in\gamma_{r_0}$ на змінну $Z\in \overline{\mathcal{U}_{r_0}}$.
	Для цього зауважимо, що для кожної фіксованої точки $S\in  \gamma_{r_0}$ функція $1/d(S,Z)$ за змінною $Z$
	є голоморфною в облаcті  $\mathcal{U}_{r_0}$ і, крім того, неперервено продовжується  на  $\gamma_{r_0}\setminus\{S\}$.
	Неперервна продовжуваність цієї функції у точку $Z=S$ випливає з існування тілесної похідної відображення $\sigma$ у цій точці і нерівності $\sigma'(S)\ne 0$.
	Отже, для кожної фіксованої точки $S\in  \gamma_{r_0}$ функція $1/d(S,Z)$ за змінною $Z$ є неперервною в $\overline{\mathcal{U}_{r_0}}$ і голоморфною в області $ \mathcal{U}_{r_0}$.
	
	Тепер, спираючись на принцип максимума для голоморфних функцій, робимо висновок  про те, що нерівність \eqref{Nyzh-oz-cont} залишаються в силі при заміні
	змінної $T\in\gamma_{r_0}$ на змінну $Z\in \overline{\mathcal{U}_{r_0}}$.
	
	Лему доведено.

	\vskip 2mm

	Через $s[Z,T]$ позначимо відрізок з початком у точці $Z$ і кінцем у точці $T$.

	\medskip
	\begin{lemma}
		\label{relatSigCoTil}
		{\it Нехай  контурна похідна $\sigma'$  неперервна та не обертається в   нуль на  $\Gamma \setminus \mathcal{X}$ і, крім того,  тілесна похідна  $\sigma'(Z)$ задовольняє оцінку
			\begin{equation}
				\label{InEqTilConDeCorP}
				\bigl|\sigma'(Z)\bigr|  \le c \prod_{j=1}^{m}\left|Z-X_j\right|^{\beta_j}\qquad \forall\, Z \in \overline{\mathcal{U}}\setminus  \mathcal{X},
			\end{equation}
			де $\beta_j \in (-1;1]$ і стала $c$ не залежить від $Z$.
			
			Тоді виконується нерівність
			\begin{equation}
				\label{sigmKIneq}
				\sum_{k=1}^2
				\bigl|d_{k}(S,Z)\bigr|
				\le c \prod_{j=1}^{m}\Bigl(\left|S-X_j \right|+\left|Z-X_j \right|\Bigr)^{\beta_j}
				\qquad \forall\, Z \in \overline{\mathcal{U}}\setminus  \mathcal{X} \quad  \forall\, S \in \Gamma \setminus \mathcal{X},
			\end{equation}
			де стала\, $c$\, не залежить  від   $S$ і $Z$.	
		}
	\end{lemma}
	
	{\bf \em Доведення.}
	Зафіксуємо $r_0\in \left(0;\mathcal{R}_{0}/4\right)$ і
	розглянемо наступні випадки розміщення точок $Z\in\overline{\mathcal{U}}\setminus\mathcal{X}$ і $S\in\Gamma\setminus\mathcal{X}$.
	Зазначимо, що скрізь у доведенні через\, $c$\, позначено сталі, значення яких не залежить  від   $S$ і $Z$,
	але, взагалі кажучи, різні навіть у межах одного ланцюжка нерівностей.
	
	{  \bf Випадок 1.}
	Нехай 
	$Z\in \overline{\mathcal{U}}_{r_0}(X_j)$ і $S \in \Gamma_{2r_0}(X_j)$ при деякому $X_j\in\mathcal{X}$
	і виконується умова $|S-Z|\le (1/2)|Z-X_j|$. Тоді $(1/2) |Z-X_j|\le |S-X_j|\le (3/2) |Z-X_j|$ і, враховуючи співвідношення \eqref{InEqTilConDeCorP}, для кожного $k\in\{1,2\}$ маємо оцінки 
	\begin{multline*}
		\bigl| \sigma_{k}\left(S\right)-
		\sigma_{k}\left(Z\right) \bigr| \le \bigl| \sigma\left(S\right)-
		\sigma\left(Z\right) \bigr| \le \int\limits_{s[Z,S]} \bigl|\sigma'\left(S\right) \bigr| |dS| \le  \\
		\le  c\, \int\limits_{s[Z,S]} |S-X_j|^{\beta_j} \, |dS| \le c \, |Z-X_j|^{\beta_j}|S-Z|,
	\end{multline*}
	з яких випливає нерівність \eqref{sigmKIneq} у випадку 1.
	
	{  \bf Випадок 2.} Нехай
	$Z\in \overline{\mathcal{U}}_{r_0}(X_j)$ і $S \in \Gamma_{2r_0}(X_j)$ при деякому $X_j\in\mathcal{X}$
	і виконується умова $|S-Z| > (1/2)|Z-X_j|$. Тоді, враховуючи співвідношення \eqref{InEqTilConDeCorP},
	для кожного $k\in\{1,2\}$ отримуємо
	\begin{multline*}
		\bigl| \sigma_{k}\left(S\right)-\sigma_{k}\left(Z\right) \bigr| \le
		\bigl| \sigma\left(S\right)- \sigma\left(Z\right) \bigr| \le
		\bigl| \sigma\left(S\right)-\sigma\left(X_j\right)\bigr|+\bigl|\sigma\left(Z\right)-\sigma\left(X_j\right) \bigr|\le\\
		\le \int\limits_{s[X_j,S]}\bigl|\sigma'\left(S\right) \bigr|\,|dS| +  \int\limits_{s[X_j,Z]}\bigl|\sigma'\left(S\right)\bigr|\,|dS| \le
		c \biggl(\,\int\limits_{s[X_j,S]}+\int\limits_{s[X_j,Z]}\biggr) |S-X_j|^{\beta_j} \, |dS|\le \\
		\le  c\left(\left|S- X_j \right|^{\beta_j + 1} + \left|Z- X_j \right|^{\beta_j+1}\right)  \le
			c \Bigl(\max\left\{|S-X_j|, |Z-X_j| \right\}\Bigr)^{\beta_j +1}.
		\end{multline*}
		
		Тепер за умови $|Z-X_j|\ge |S-X_j|$ використовуємо нерівність $|S-Z| > (1/2)|Z-X_j|$ і отримуємо оцінку
		\[ | d_{k}(S,Z)|\le c \,|Z-X_j|^{\beta_j }, \]
		а за умови $|Z-X_j|<|S-X_j|$ використовуємо нерівність $|S-Z|\ge (1/2)|S-X_j|$ і маємо
		\[ | d_{k}(S,Z)|\le c \,|S-X_j|^{\beta_j }, \]
		звідки випливає нерівність \eqref{sigmKIneq} у випадку 2.
		
		{  \bf Випадок 3.}
		Нехай $Z\in \overline{\mathcal{U}}_{r_0}\left(X_j\right)$ і $S \in \Gamma \setminus  \Gamma_{2r_0}\left(X_j\right)$  при деякому $X_j\in\mathcal{X}$.
		Тоді $|S-Z|\ge r_0$  і для кожного $k\in\{1,2\}$ виконуються нерівності	
		\[\left| d_{k}\left(S,Z\right)\right| \le \frac{\left| \sigma(S)-\sigma(Z)\right|}{|S-Z|} \le \frac{2\,\max\limits_{S\in\overline{\mathcal{U}}} |\sigma(S)|}{r_0}\le c\, |S-X_j|^{\beta_j}, \]
		наслідком яких є нерівність \eqref{sigmKIneq} у цьому випадку.	
		
    	{  \bf Випадок 4.}
		Нехай  $Z\in \overline{\mathcal{U}} \setminus \bigcup\limits_{j=1}^{m}\,
		\overline{\mathcal{U}}_{r_0}\left(X_j\right)$ і $S \in\Gamma_{r_0/2}\left(X_k\right)$  при деякому $X_j\in\mathcal{X}$. Тоді
		$|S-Z|\ge r_0/2$  і для кожного $k\in\{1,2\}$ виконуються нерівності	
		\[\left| d_{k}\left(S,Z\right)\right| \le \frac{\left| \sigma(S)-\sigma(Z)\right|}{|S-Z|} \le
		\frac{2\,\max\limits_{S\in\overline{\mathcal{U}}} |\sigma(S)|}{r_0/2}\le c\, |Z-X_j|^{\beta_j}, \]
		наслідком яких є нерівність \eqref{sigmKIneq} у цьому випадку.

		{ \bf Випадок 5.} 	Нехай $Z\in \overline{\mathcal{U}} \setminus\,  \bigcup\limits_{j=1}^{m}\,
		\overline{\mathcal{U}}_{r_0}\left(X_j\right)$ і $S \in \Gamma \setminus \,
		\bigcup\limits_{j=1}^{m}\,    \Gamma_{r_0/2}\left(X_j\right)$.
		Тоді в силу леми~\ref{DbiggerZe} виконується оцінка вигляду \eqref{dTypPokhInUbezX}, наслідком якої
		є нерівність \eqref{sigmKIneq} у цьому випадку.
		
		Лему доведено.

		\vskip 2mm
		
		\medskip
		\begin{lemma}
			\label{relatSigCoTLow}
			{\it Нехай  контурна похідна $\sigma'$  неперервна та не обертається в  нуль на  $\Gamma \setminus \mathcal{X}$.
				Якщо при деякому $r_0\in \left(0, \mathcal{R}_{0}/4\right)$ і кожному $X_j\in\mathcal{X}$, $j=\overline{1,m}$,
				функція $d(S,Z)$ задовольняє оцінку
				\begin{equation} \label{IneQLow}
					\bigl|d(S,Z)\bigr|
					\ge c \Bigl(\left|S-X_j \right|+\left|Z-X_j \right|\Bigr)^{\beta_j}\qquad \forall\, Z \in
					\overline{\mathcal{U}}_{r_0}\left(X_j\right)\setminus \mathcal{X}
					\quad \forall\, S \in  \Gamma_{2r_0}
					\left(X_j\right)\setminus \mathcal{X},
				\end{equation}
				де $\beta_j \in (-1;1]$ і стала $c$ не залежить від $S$ і $Z$,
				то виконується також нерівність
				\begin{equation} \label{d-oz-znyzu-zag}
					\bigl|d(S,Z)\bigr| \ge
					c \prod_{j=1}^{m}\Bigl(\left|S-X_j \right|+\left|Z-X_j \right|\Bigr)^{\beta_j}
					\qquad \forall\, Z \in \overline{\mathcal{U}}\setminus  \mathcal{X} \quad  \forall\, S \in \Gamma \setminus \mathcal{X},
				\end{equation}
				де стала $c$ не залежить  від   $S$ і $Z$.	}
		\end{lemma}
		
		{\bf \em Доведення.} Для доведення нерівності \eqref{d-oz-znyzu-zag} залишається розглянути випадки 3 --- 5  розміщення точок $Z\in\overline{\mathcal{U}}\setminus\mathcal{X}$ і $S\in\Gamma\setminus\mathcal{X}$ з доведення леми \ref{relatSigCoTil}, оскільки у розглянутих там випадках 1 і 2 за припущенням леми виконується оцінка \eqref{IneQLow}, з якої випливає нерівність \eqref{d-oz-znyzu-zag} у цих випадках.
		
		У випадку 5 в силу леми~\ref{Dbig-oz-znyzu} виконується оцінка вигляду \eqref{d-oz-znyzu}, наслідком якої
		є нерівність \eqref{d-oz-znyzu-zag} у цьому випадку.
		
		У випадку 3 виконується нерівність $|S-Z|>r_0$, а у випадку 4 --- нерівність $|S-Z|>r_0/2$. Тому в цих випадках
		подібно до оцінки \eqref{Nyzh-oz-homeomorph} приходимо до висновку, що існує додатна стала $c_0(r_0)$, залежна від $r_0$, така, що виконується нерівність
		\[ |d(S,Z)|\ge c_0(r_0),\]
		наслідком якої є нерівність \eqref{d-oz-znyzu-zag} у випадках 3 і 4.
		
		Лему доведено.
		
		Для пари різних точок $S,T\in \Gamma$ через $\Gamma[S,T]$ позначимо ту дугу кривої $\Gamma$ з початком у  $S$ і кінцем у $T$,
		яка має небільшу довжину $l\left(\Gamma[S,T]\right)$.
				
		Очевидно, що виконується подвійна нерівність
		\begin{equation}
			\label{duga-khorda-kolo}
			1 \le \frac{l\left(\Gamma[S,T]\right)}{|S-T|}\le \frac{\pi}{2} \qquad \forall\,S,T\in \Gamma : S\ne T.
		\end{equation}

		Функція $\omega \colon (0; +\infty) \longrightarrow (0; +\infty)$ називається {\em напівадитивною}, якщо
		\[\omega(\eta_1+\eta_2)\le \omega(\eta_1)+\omega(\eta_2) \qquad \forall\,\eta_1>0\quad \forall\,\eta_2>0.\]
		Відомо (див., наприклад, \cite[с. 175]{Dz-Shev-6-3}), що для неспадної напівадитивної функції $\omega$ виконується нерівність
		\begin{equation}
			\label{napivad-1}
			\omega(\lambda\eta)\le (\lambda+1) \omega(\eta) \qquad \forall\,\lambda>0\quad \forall\,\eta>0,
		\end{equation}
		яка уточнюється у випадку натурального $\lambda=n$:
		\begin{equation}
			\label{napivad-2}
			\omega(n\eta)\le n \omega(\eta) \qquad  \forall\,\eta>0.
		\end{equation}
		
		\medskip
		\begin{lemma}
			\label{relatSigCoTLow2}
			{\it Нехай  контурна похідна $\sigma'$  неперервна та не обертається в   нуль на  $\Gamma \setminus \mathcal{X}$, а  тілесна похідна $\sigma'(Z)$ задовольняє оцінку
				\eqref{InEqTilConDeCorP}. Нехай, крім того,  
				при деякому $r_0 \in \left(0;\mathcal{R}_0/4\right)$ виконуються оцінки	
				\begin{multline}
					\label{ots-SiDerLoc1}
					\left|\sigma'\left(Z\right)-
					\sigma'\left(Z_0\right)\right| 
					\le c\,
					\frac{\omega_1^{\ast}\left(|Z-Z_0|\right)}{\omega_1^{\ast}\left(\frac{1}{2} \mathcal{R}(Z_0)\right)}
					\prod_{j=1}^{m} \left| Z-X_j \right|^{\beta_j}\\
					\qquad  \forall Z_0 \in
					\Gamma_{r_0} \left\langle \mathcal{X}\right\rangle \setminus \mathcal{X}
					\quad
					\forall Z \in  \Big(\Gamma\cup s[0,Z_0]\Big)_{ \mathcal{R}\left(Z_0\right)/2}\left(Z_0\right),
				\end{multline}
				\begin{multline}
					\label{ots-SiDerLoc}
					\left|\sigma'\left(Z\right)-
					\sigma'\left(Z_0\right)\right| 
					\le c\,
					\frac{\omega_1\left(|Z-Z_0|\right)}{\omega_1\left(\frac{1}{2} \mathcal{R}(Z_0)\right)}
					\prod_{j=1}^{m} \left| Z-X_j \right|^{\beta_j}\\
					\qquad  \forall Z_0 \in \Gamma \setminus
					\Gamma_{r_0} \left\langle \mathcal{X}\right\rangle
					\,\,\,
					\forall Z \in    \Big(\Gamma\cup s[0,Z_0]\Big)_{ \mathcal{R}\left(Z_0\right)/2}\left(Z_0\right),
				\end{multline}
				в яких $\beta_j$, $j=\overline{1,m}$, ті ж самі числа, що й в оцінці \eqref{InEqTilConDeCorP}, стала
				$c$ не залежить від  $Z$ і $Z_0$, $\omega_1 \colon (0; +\infty) \longrightarrow (0; +\infty)$ і $\omega_1^{\ast} \colon (0; +\infty) \longrightarrow (0; +\infty)$~---
				неспадні обмежені напівадитивні  функції такі, що  $\omega_1$ задовольняє умову \eqref{dini}, а  $\omega_1^{\ast}$~--- умову \eqref{dini-om1Pr}.
				
				Тоді  виконується оцінка
				\begin{multline}
					\label{Lemm8MainOts}
					\sum_{k=1}^{2}\Bigl|\sigma_{k}'(S)- d_{k}(S,Z) \Bigr| \le c \,\min \left\{\frac{\omega\left(|S-Z|\right)}{\omega\left(\frac{1}{2} \mathcal{R}(Z)\right)}, 1\right\}\times\\
					\times \left(\prod_{j=1}^{m}\left|S-X_j \right|^{\beta_j}+
					\prod_{j=1}^{m} \Bigl(\left|S-X_j\right|+
					\left|Z-X_j\right|\Bigr)^{\beta_j} \right)\qquad \forall Z\in \overline{\mathcal{U}}\setminus \mathcal{X}\quad \forall S\in \Gamma \setminus \mathcal{X},
				\end{multline}
				де
				\begin{equation}
					\label{DefOm-tsetrZ}
					\omega\left(\varepsilon\right):=
					\begin{cases}
						\omega_{1}^{\ast}\left(\varepsilon\right),  &\text{якщо $\mathcal{R}(Z) \le {r_0}/{2}$,}\\
						\omega_{1}\left(\varepsilon\right) +
						\omega_{1}^{\ast}\left(\varepsilon\right),  &\text{якщо $\mathcal{R}(Z) > {r_0}/{2}$,}
					\end{cases}
				\end{equation}
				а стала $c$ не залежить від $S$ і $Z$.
			}
		\end{lemma}
		
		{\bf \em Доведення.} Доведемо спочатку оцінку \eqref{Lemm8MainOts} для $Z=Z_0\in\Gamma \setminus \mathcal{X}$.
		
		У випадку, коли $S\in \Gamma$ і  $|S-Z_0|<\mathcal{R}(Z_0)/2$,
		для будь-якого $k\in\{1,2\}$ і $T\in \Gamma_{\mathcal{R}(Z_0)/2}(Z_0)$, враховуючи нерівність \eqref{duga-khorda-kolo},
		отримуємо співвідношення
		\begin{multline*}
			\Big|d_{k}\left(S,T\right)-
			d_{k}\left(S,Z_0\right)\Big| =\\
			=\left|\frac{1}{S-T} \int\limits_{\Gamma\left[T,S\right]}\Big(\sigma_{k}'(S)- \sigma_{k}'(Z_0)\Big)  \, dS-
			\frac{1}{S-Z_0} \int\limits_{\Gamma\left[Z_0, S\right]}\Big(\sigma_{k}'(S)- \sigma_{k}'(Z_0)\Big)  \, dS \right| \le \\
			\le \left| \left(\frac{1}{S-T}-\frac{1}{S-Z_0}\right)\int\limits_{\Gamma\left[T, S\right]}\Big(\sigma_{k}'(S)- \sigma_{k}'(Z_0)\Big)\,dS\right|+\left|\frac{1}{S-Z_0}\int\limits_{\Gamma\left[T, Z_0\right]}\Big(\sigma_{k}'(S)- \sigma_{k}'(Z_0)\Big)  \, dS\right|\le \\
			\le
			\frac{|T-Z_0|}{|S-Z_0|}\frac{l\left(\Gamma[S,T]\right)}{|S-T|}\sup_{S_{\ast}\in \Gamma[S,T]}\Bigl|\sigma_{k}'\left(S_{\ast}\right)-\sigma_{k}'\left(Z_0\right)\Bigr|+\frac{l\left(\Gamma[T,Z_0]\right)}{|S-Z_0|}\sup_{S_{\ast}\in \Gamma[T,Z_0]}\Bigl|\sigma_{k}'\left(S_{\ast}\right)-\sigma_{k}'\left(Z_0\right)\Bigr| \le \\
			\le \frac{|T-Z_0|}{|S-Z_0|}
			\left(\frac{l\left(\Gamma[S,T]\right)}{|S-T|}\sup_{S_{\ast}\in \Gamma[S,T]} \Bigl|\sigma_{k}'\left(S_{\ast}\right)-\sigma_{k}'\left(Z_0\right)\Bigr|+  
			\frac{l\left(\Gamma[T,Z_0]\right)}{|T-Z_0|} \sup_{S_{\ast}\in \Gamma[T,Z_0]} \Bigl|\sigma_{k}'\left(S_{\ast}\right)-\sigma_{k}'\left(Z_0\right)\Bigr|\right)\le
		\end{multline*}
		\begin{equation}
			\label{supdkLoc}
			\le \frac{\pi}{2} \, \frac{\left| T-Z_0\right|}{\left| S-Z_0\right|} \left(\sup_{S_{\ast}\in \Gamma[S,T]} \Bigl|\sigma_{k}'\left(S_{\ast}\right)-\sigma_{k}'\left(Z_0\right)\Bigr|+ \sup_{S_{\ast}\in \Gamma[T,Z_0]} \Bigl|\sigma_{k}'\left(S_{\ast}\right)-\sigma_{k}'\left(Z_0\right)\Bigr|\right).
		\end{equation}
		Переходячи до границі в нерівності \eqref{supdkLoc} при $T \to S$, 
		отримуємо нерівність
		\[\left|\sigma_{k}'\left(S\right)-
		d_{k}\left(S,Z_0\right)\right| \le \pi \sup_{S_{\ast}\in \Gamma[S,Z_0]}\Bigl|\sigma_{k}'\left(S_{\ast}\right)-\sigma_{k}'\left(Z_0
		\right)\Bigr|,\]
		і тепер у розглянутому випадку оцінка \eqref{Lemm8MainOts}
		є наслідком оцінок \eqref{ots-SiDerLoc1} і \eqref{ots-SiDerLoc}.
		
		У випадку, коли $S\in \Gamma\setminus \mathcal{X}$ і  $|S-Z_0|\ge\mathcal{R}(Z_0)/2$,
		для кожного $k\in\{1,2\}$ виконується нерівність
		\[ \left|\sigma_{k}'\left(S\right)-
		d_{k}\left(S,Z_0\right)\right| \le \left|\sigma'\left(S\right) \right|+\left|d_{k}\left(S,Z_0\right) \right|\]
		і оцінка \eqref{Lemm8MainOts}
		є наслідком оцінок	\eqref{InEqTilConDeCorP} і \eqref{sigmKIneq}.

		Доведемо тепер оцінку \eqref{Lemm8MainOts} для $Z\in \mathcal{U}$.
		
		Нехай $Z_{\ast}$ --- точка кола $\Gamma$, найближча до точки $Z$.
		
		При $|Z-Z_*|<\mathcal{R}(Z)/8$ розглянемо два випадки розміщення точки $S\in \Gamma \setminus \mathcal{X}$.
		
		У випадку, коли $S\in \Gamma$ і  $|S-Z|<\mathcal{R}(Z)/4$, для кожного $k\in\{1,2\}$ використовуємо нерівність
		\begin{equation}	\label{proizv-dkLoc}
			\bigl| \sigma_{k}'\left(S\right)-d_{k}\left(S,Z\right)\bigr|\le \bigl|  \sigma_{k}'\left(S\right)-d_{k}\left(S,Z_{\ast}\right)\bigr| +
			\bigl| d_{k}\left(S, Z_{\ast}\right)-d_{k}\left(S,Z\right)\bigr|.
		\end{equation}
		
		Оскільки для
		всіх $S\in \Gamma : |S-Z|<\mathcal{R}(Z)/4$ виконуються співвідношення
		\[|S-Z_*|\le |S-Z|+|Z-Z_*|< \frac{1}{4}\,\mathcal{R}(Z)+\frac{1}{8}\,\mathcal{R}(Z)<\frac{3}{8}\cdot\frac{8}{7}\,\mathcal{R}(Z_*)<\frac{1}{2}\,\mathcal{R}(Z_*),\]
		то для першого доданка з правої частини нерівності \eqref{proizv-dkLoc} 
		вище доведена оцінка вигляду \eqref{Lemm8MainOts}, в якій замість точки $Z=Z_0\in\Gamma \setminus \mathcal{X}$ слід підставити $Z_*$.
		
		При оцінюванні другого доданка з правої частини нерівності \eqref{proizv-dkLoc} з урахуванням нерівності \eqref{duga-khorda-kolo}
		отримуємо співвідношення
		\begin{multline*}
			\bigl| d_{k}\left(S, Z_{\ast}\right)-d_{k}\left(S,Z\right)\bigr|=\\
			= \left|\frac{1}{S-Z_{\ast}} \int\limits_{\Gamma\left[Z_{\ast},S\right]} \Big(\sigma_{k}'(S)- \sigma_{k}'(Z_*)\Big)\,dS-
			\frac{1}{S-Z} \int\limits_{\Gamma\left[Z_{\ast}, S\right]\bigcup s\left[Z,Z_{\ast}\right]} \Big(\sigma_{k}'(S)- \sigma_{k}'(Z_*)\Big)\,dS \right| \le \\
			\le \left| \left(\frac{1}{S-Z_*}-\frac{1}{S-Z}\right)\int\limits_{\Gamma\left[Z_*, S\right]}\Big(\sigma_{k}'(S)- \sigma_{k}'(Z_*)\Big)\,dS\right|+\left|\frac{1}{S-Z}\int\limits_{s\left[Z,Z_{\ast}\right]}\Big(\sigma_{k}'(S)- \sigma_{k}'(Z_*)\Big)  \, dS\right|\le 
		\end{multline*}
		\[\le \sup_{S_{\ast}\in \Gamma\left[Z_{\ast}, S\right]\bigcup s\left[Z,Z_{\ast}\right]}
		\Bigl|\sigma_{k}'\left(S_{\ast}\right)-\sigma_{k}'\left(Z_*\right)\Bigr| \frac{|Z-Z_*|}{|S-Z|}
		\left(\frac{l\left(\Gamma[Z_*,S]\right)}{|S-Z_*|}+1  \right) \le \]
		\begin{equation}
			\label{otsdkLonNearZast}
			\le
			\frac{|Z-Z_*|}{|S-Z|}
			\left(\frac{\pi }{2}+1  \right) \sup_{S_{\ast}\in \Gamma\left[Z_{\ast}, S\right]\bigcup s\left[Z,Z_{\ast}\right]}
			\Bigl|\sigma_{k}'\left(S_{\ast}\right)-\sigma_{k}'\left(Z_*\right)\Bigr|,
		\end{equation}
		наслідком яких є нерівність
		\[\bigl| d_{k}\left(S, Z_{\ast}\right)-d_{k}\left(S,Z\right)\bigr|
		\le \Big(\frac{\pi}{2}+1\Big)\sup_{S_{\ast}\in \Gamma\left[Z_{\ast}, S\right]\bigcup s\left[Z,Z_{\ast}\right]} \Bigl|\sigma_{k}'\left(S_{\ast}\right)-\sigma_{k}'\left(Z_*\right)\Bigr|. \]
		Отже, оцінка другого доданка з правої частини нерівності \eqref{proizv-dkLoc} має такий же вигляд, як і оцінка першого доданка з правої частини цієї нерівності.
		
		Покажемо, що у вказаних оцінках доданків з правої частини нерівності \eqref{proizv-dkLoc} можна замінити точку $Z_*$ точкою $Z$.
		Для цього зазначимо, що для всіх $S_{\ast} \in \Gamma\left[Z_{\ast},S\right]$ виконуються нерівності
		$\left|S_{\ast}-Z_{\ast} \right| \le \left|S- Z_{\ast}\right|\le 2|S-Z|$, а для всіх $S_{\ast} \in s\left[Z, Z_{\ast}\right]$ --- нерівності
		$\left|S_{\ast}-Z_{\ast} \right| \le |Z_{\ast}-Z| \le |S-Z|$.
		Оскільки також $\mathcal{R}(Z_*)>(7/8)\mathcal{R}(Z)$, то, враховуючи нерівності \eqref{napivad-1} і \eqref{napivad-2},
		при всіх $S_{\ast}\in \Gamma\left[Z_{\ast}, S\right]\bigcup s\left[Z,Z_{\ast}\right]$  отримуємо
		\[\frac{\omega\left(\left|S_{\ast}-Z_{\ast}\right|\right)}{\omega\left(\frac{1}{2}\mathcal{R}(Z_{\ast})\right)}\le
		\frac{\omega\left(2|S-Z|\right)}{\omega\left(\frac{7}{16}\mathcal{R}(Z)\right)}=
		\frac{\omega\left(2|S-Z|\right)\,\omega\left(\frac{8}{7}\cdot\frac{7}{16}\mathcal{R}(Z)\right)}{\omega\left(\frac{1}{2}\mathcal{R}(Z)\right)\,
			\omega\left(\frac{7}{16}\mathcal{R}(Z)\right)}\le\]
		\[\le 2\left(\frac{8}{7}+1\right)\frac{\omega\left(|S-Z|\right)}{\omega\left(\frac{1}{2}\mathcal{R}(Z)\right)}=
		2\left(\frac{8}{7}+1\right)\,\min \left\{\frac{\omega\left(|S-Z|\right)}{\omega\left(\frac{1}{2} \mathcal{R}(Z)\right)}, 1\right\}. \]
		Тепер можна стверджувати, що у випадку, коли $S\in \Gamma$ і  $|S-Z|<\mathcal{R}(Z)/4$, виконується оцінка \eqref{Lemm8MainOts}.

		У випадку, коли $S\in \Gamma\setminus \mathcal{X}$ і $|S-Z|\ge\mathcal{R}(Z)/4$,
		для кожного $k\in\{1,2\}$ скористаємось нерівністю
		\[ \left|\sigma_{k}'\left(S\right)-
		d_{k}\left(S,Z\right)\right| \le \left|\sigma'\left(S\right) \right|+\left|d_{k}\left(S,Z\right) \right|\]
		і оцінками	\eqref{InEqTilConDeCorP} і \eqref{sigmKIneq}.
		Крім того, з урахуванням нерівності \eqref{napivad-2} отримаємо оцінку
		\[\frac{\omega\left(|S-Z|\right)}{\omega\left(\frac{1}{2} \mathcal{R}(Z)\right)}\ge
		\frac{\omega\left(\frac{1}{4} \mathcal{R}(Z)\right)}{\omega\left(2\cdot\frac{1}{4} \mathcal{R}(Z)\right)}\ge \frac{1}{2}, \]
		з якої випливає нерівність
		\begin{equation}	\label{oz-vidn-mod-nepr-znyzu}
			1\le 2\,\min \left\{\frac{\omega\left(|S-Z|\right)}{\omega\left(\frac{1}{2} \mathcal{R}(Z)\right)}, 1\right\}.
		\end{equation}
		Отже, в розглянутому випадку оцінка \eqref{Lemm8MainOts} також виконується.
		
		Нарешті, при $|Z-Z_*|\ge\mathcal{R}(Z)/8$,
		враховуючи нерівність \eqref{napivad-2}, отримуємо оцінку
		\[\frac{\omega\left(|S-Z|\right)}{\omega\left(\frac{1}{2} \mathcal{R}(Z)\right)}\ge
		\frac{\omega\left(\frac{1}{8} \mathcal{R}(Z)\right)}{\omega\left(4\cdot\frac{1}{8} \mathcal{R}(Z)\right)}\ge \frac{1}{4},\]
		і оцінка \eqref{Lemm8MainOts} встановлюється так, як і в попередньому випадку.
		Лему доведено.
		
		\vskip 2mm

		\begin{lemma}
			\label{DiffDkNearZ0}
			{\it За умов  леми~\ref{relatSigCoTLow2}  виконується нерівність
				\begin{equation}
					\label{IneqdifDZ0}
					\begin{split}
						\sum\limits_{k=1}^{2} \left|d_{k} \left(S,Z\right)- d_{k} \left(S,Z_0\right) \right| \le  c \, \frac{\left|Z-Z_0\right|}{\left|S-Z_0\right|} \,\min \left\{\frac{\omega\left(|S-Z_0|\right)}{\omega\left(\frac{1}{2} \mathcal{R}(Z_0)\right)}, 1\right\}\times \\
						\times \left(\,\prod_{j=1}^m \Big(\left|S - X_j \right| + \left|Z_0 - X_j \right| \Big)^{\beta_j}
						+  \prod_{j=1}^m \Big(\left|Z - X_j \right| + \left|Z_0 - X_j \right| \Big)^{\beta_j} \right)   \\			
						\forall Z_0 \in \Gamma \setminus \mathcal{X} \quad \forall  Z \in \overline{\mathcal{U}}\setminus  \mathcal{X}  \quad
						\forall S \in \Gamma \setminus \mathcal{X}: |S-Z_0| \ge 2|Z-Z_0|,
					\end{split}
				\end{equation}
				де функція $\omega$ визначена рівністю \eqref{DefOm-tsetrZ},
				а стала $c$ не залежить від $S$, $Z_0$ і $Z$.	}
		\end{lemma}
		
		{\bf \em Доведення.}
		Для довільного $Z_0\in\Gamma \setminus \mathcal{X}$ розглянемо спочатку випадок, коли
		\[	Z\in \left(\overline{\mathcal{U}}_{\mathcal{R}\left(Z_0\right)/16}\right)\left(Z_0\right) \setminus \mathcal{X},\,\, S\in \Gamma \setminus \mathcal{X}:	2\left|Z-Z_0\right|\le \left|S-Z_0\right| \le \left(1/8\right)\mathcal{R}\left(Z_0\right). \]

		Нехай $Z_{\ast}$ --- точка кола $\Gamma$, найближча до точки $Z$.
		
		Для кожного $k=1,2$ використаємо нерівність
		\begin{equation}	\label{dkLocZast}
			\bigl| d_{k}\left(S,Z\right)-d_{k}\left(S,Z_0\right)\bigr|\le \bigl|  d_{k}\left(S,Z\right)-d_{k}\left(S,Z_{\ast}\right)\bigr| +
			\bigl| d_{k}\left(S, Z_{\ast}\right)-d_{k}\left(S,Z_0\right)\bigr|.
		\end{equation}
		
		Враховуючи нерівність $\mathcal{R}\left(Z_0\right) \le (16/15)\mathcal{R}\left(Z\right)$, отримуємо співвідношення
		\[ \left|Z-Z_{\ast}\right|\le \left|Z-Z_0\right| \le  (1/16) \mathcal{R}\left(Z_0\right) \le (1/15)\mathcal{R}\left(Z\right)<(1/8)\mathcal{R}\left(Z\right),\]
		\[|S-Z|\le |S-Z_0|+|Z-Z_0|\le (1/8) \mathcal{R}(Z_0)+(1/16) \mathcal{R}(Z_0) 
		\le (1/5) \mathcal{R}(Z)<(1/4) \mathcal{R}(Z). \]
		Тому для першого доданка з правої частини нерівності \eqref{dkLocZast} виконується оцінка
		\eqref{otsdkLonNearZast}.
		Для другого доданка з правої частини нерівності  \eqref{dkLocZast} використовуємо оцінку
		\eqref{supdkLoc}, в якій замість точки  $T$ слід підставити $Z_{\ast}$,
		оскільки при цьому 
		\[\left|S-Z_0\right|\le (1/8)\mathcal{R}\left(Z_0\right),\quad
		\left|Z_{\ast}-Z_0\right| \le
		2\left|Z-Z_0\right| \le (1/8)\mathcal{R}\left(Z_0\right).\]
		
		В результаті маємо
		\[ \bigl| d_{k}\left(S,Z\right)-d_{k}\left(S,Z_0\right)\bigr|\le \frac{|Z-Z_*|}{|S-Z|}
		\left(\frac{\pi}{2}+1  \right)  \sup_{S_{\ast}\in \Gamma\left[Z_{\ast}, S\right]\bigcup s\left[Z,Z_{\ast}\right]}
		\Bigl|\sigma_{k}'\left(S_{\ast}\right)-\sigma_{k}'\left(Z_*\right)\Bigr| + \]
		\[+\frac{\pi}{2} \, \frac{\left| Z_{\ast}-Z_0\right|}{\left| S-Z_0\right|} \left(\sup_{S_{\ast}\in \Gamma[S,Z_{\ast}]} \Bigl|\sigma_{k}'\left(S_{\ast}\right)-\sigma_{k}'\left(Z_0\right)\Bigr|+ \sup_{S_{\ast}\in \Gamma[Z_{\ast},Z_0]} \Bigl|\sigma_{k}'\left(S_{\ast}\right)-\sigma_{k}'\left(Z_0\right)\Bigr|\right).\]
		
		Тепер у розглянутому випадку оцінка \eqref{IneqdifDZ0} є наслідком оцінок \eqref{ots-SiDerLoc1} і \eqref{ots-SiDerLoc}, співвідношень
		\[\left|S-Z\right|\ge \left|S-Z_0\right|/2, \quad  \left|Z-Z_{\ast}\right|\le \left|Z-Z_0\right|, \quad
		\left|Z_{\ast}-Z_0\right| \le2\left|Z-Z_0\right|, \quad \mathcal{R}\left(Z_{\ast}\right)\ge (7/8)\mathcal{R}\left(Z_{0}\right),\]
		\[\left|S_{\ast}-Z_{\ast}\right| \le 2\left|S-Z_0\right| \qquad \forall\, S_{\ast}\in  \Gamma[S,Z_{\ast}]\cup s[Z_{\ast},Z] \]
		і властивостей монотонності та напівадитивності функцій $\omega_1$ і $\omega_1^*$.
		
		Нарешті, для довільного $Z_0\in\Gamma \setminus \mathcal{X}$ розглянемо випадки, коли
		\[	Z\in \left(\overline{\mathcal{U}}_{\mathcal{R}\left(Z_0\right)/16}\right)\left(Z_0\right) \setminus \mathcal{X},\,\, S\in \Gamma \setminus \mathcal{X}:
		\left|S-Z_0\right| > \left(1/8\right)\mathcal{R}\left(Z_0\right)\]
		або
		\[		Z\in \overline{\mathcal{U}} \setminus  \left(\overline{\mathcal{U}}_{\mathcal{R}\left(Z_0\right)/16}\right)\left(Z_0\right)\setminus  \mathcal{X}; \,\, S\in \Gamma \setminus \mathcal{X}:	\left|S-Z_0\right|\ge 2 \left|Z-Z_0\right|.	\]

		Враховуючи нерівність  $|S-Z| \ge |S-Z_{0}|/2$, що виконується в цих випадках, і оцінку \eqref{sigmKIneq}, для кожного  $k\in \{1,2\}$  отримуємо   співвідношення
		\[\bigl| d_{k}\left(S, Z\right)-d_{k}\left(S,Z_0\right)\bigr|=\left|\frac{(Z-Z_0)\big(\sigma_{k}(S)-\sigma_{k}(Z_0)\big)}{(S-Z)(S-Z_0)}+
		\frac{\sigma_{k}(Z_0)-\sigma_{k}(Z)}{S-Z} \right|=\]
		\[= \frac{|Z-Z_0|}{|S-Z|}\, \Bigl|d_{k}\left(S,Z_0\right) - d_{k}\left(Z_0,Z\right)\Bigr|
		\le \frac{2\,|Z-Z_0|}{|S-Z_0|}\,\Bigl(\left|d_{k}\left(S,Z_0\right)\right|+\left|d_{k}\left(Z_0,Z\right)\right|\Bigr)\le\]
		\begin{equation}\label{spiv-dlia-rozd-rizn}
			\le \frac{2\,|Z-Z_0|}{|S-Z_0|}\, \left(\,\prod_{j=1}^m \Big(\left|S - X_j \right| + \left|Z_0 - X_j \right| \Big)^{\beta_j}
			+  \prod_{j=1}^m \Big(\left|Z - X_j \right| + \left|Z_0 - X_j \right| \Big)^{\beta_j} \right).
		\end{equation}
		
		З урахуванням того, що у зазначених випадках $|S-Z_0| > (1/8)\mathcal{R}(Z_0)$,
		подібно до нерівності \eqref{oz-vidn-mod-nepr-znyzu} отримуємо нерівність
		\begin{equation} \label{ge1fr8}
			1\le 4\,\min \left\{\frac{\omega\left(|S-Z_0|\right)}{\omega\left(\frac{1}{2} \mathcal{R}(Z_0)\right)}, 1\right\},
		\end{equation}
		і тепер оцінка \eqref{IneqdifDZ0} випливає з нерівностей \eqref{spiv-dlia-rozd-rizn} і \eqref{ge1fr8}.
		Лему доведено.
		
%		\section{Зведення кусково-неперервної (1-3) задачі до системи інтегральних рівнянь}
{\bf 4. Зведення кусково-неперервної (1-3) задачі до системи інтегральних рівнянь.}
				Очевидно, що функція
		\begin{equation}\label{widetildSigm}
			\widetilde{\sigma}(Z):=   \sigma_1(Z)\, e_1 +  \sigma_2(Z) \, e_2 \qquad \forall\, Z \in \overline{\mathcal{U}}
		\end{equation}
		визначає гомеоморфне відображення  $\overline{\mathcal{U}}$ на $\overline{D_{\zeta}}$.
		Як і відображення  $\sigma$, відображення (\ref{widetildSigm}) також має неперервну контурну похідну на множині
		$\Gamma \setminus  \mathcal{X}$.
		Враховуючи рівності \eqref{tabl-umn}, \eqref{rho}, отримуємо
		\[	\widetilde{\sigma}'\left(S\right) =
			\sigma'\left(S\right) - \frac{i}{2}\,\sigma_{2}'\left(S\right)\rho \qquad
			\forall\, S \in \Gamma \setminus  \mathcal{X}.\]

		Для $\tau = \widetilde{\sigma}(S)\in \partial D_{\zeta} \setminus \Upsilon_{\zeta}$ і
		$\zeta=\widetilde{\sigma}(Z)\in D_{\zeta}$ справедлива рівність
		\[	(\tau-\zeta)^{-1}=\frac{1}{\sigma(S)-\sigma(Z)} +\frac{i\rho}{2}\frac{\sigma_2(S)-\sigma_{2}(Z)}{\left(\sigma(S)-\sigma(Z)\right)^{2}}\,.\]

		Так, як і при доведенні Леми~3  роботи \cite{GrPlHBVMbF-Flaut21} у випадку гладкого контура  $\partial D_{z}$, перетворюємо рівність  \eqref{cl-sol-cone} до вигляду
		\begin{multline} \label{BihInt13ToSumOmeg}
			\Phi(\zeta) = \frac{1}{2\pi i}\int\limits_{\Gamma}
			\frac{\widetilde{\varphi}\left(S\right)\widetilde{\sigma}'\left(S\right)}
			{\sigma\left(S\right)-\sigma\left(Z\right)}\, d S
			=\frac{e_1}{2\pi i}\int_{\Gamma}
			\frac{\Omega_{1}^{\widetilde{\varphi}}\left(S,Z\right)}{S-Z}\, d S +
			\frac{2ie_1-2e_2}{2\pi i}\int_{\Gamma}
			\frac{\Omega_{2}^{\widetilde{\varphi}}\left(S,Z\right)}{S-Z}\, d S  + \\
			+
			\frac{e_1}{4 \pi i}\int\limits_{\Gamma} \frac{\widetilde{\varphi}\left(S\right)}{S}\left(1+
			\frac{S+Z}{S-Z}\right)\, dS
			\qquad \forall\, Z \in \mathcal{U},
		\end{multline} 
		де
		\begin{equation}\label{Om1-int}
			\Omega_{1}^{\widetilde{\varphi}}\left(S,Z\right) =
			\widetilde{\varphi}\left(S\right)
			\left(\frac{\sigma'\left(S\right)}{d(S,Z)}-1\right)
			\qquad \forall\, S \in \Gamma \setminus \mathcal{X} \quad \forall\,  Z \in {\mathcal{U}},
		\end{equation}
		\begin{equation}\label{Om2-int}
			\Omega_{2}^{\widetilde{\varphi}}\left(S,Z\right) =
			\frac{\widetilde{\varphi}\left(S\right)}{2}\,
			\left(\frac{\sigma'(S) \, d_{2}\left(S,Z\right)}
			{\left(d\left(S,Z\right)\right)^{2}} -
			\frac{\sigma_{2}'(S)}{d\left(S,Z\right)}\right) \qquad
			\forall\, S \in \Gamma
			\setminus \mathcal{X} \quad \forall\,  Z \in {\mathcal{U}},
		\end{equation}
		\begin{equation} \label{VarPh}
			\widetilde{\varphi}\left(S\right) = \varphi_{1} \left(\widetilde{\sigma}\left(S\right)\right)e_1 +
			\varphi_{3} \left(\widetilde{\sigma}\left(S\right)\right)e_2 \qquad \forall\, S   \in \Gamma
			\setminus \mathcal{X} .
		\end{equation}

		Очевидно, що функції \eqref{Om1-int} та  \eqref{Om2-int} неперервно продовжуються у точки $Z\in  \Gamma
		\setminus \mathcal{X}$, $Z\ne S$. Крім того, $\Omega_{1}^{\widetilde{\varphi}}(S,Z)\to 0$ і $\Omega_{2}^{\widetilde{\varphi}}(S,Z)\to 0$  при $Z\to S\in\Gamma
		\setminus \mathcal{X}$. Тому надалі вважаємо, що дані функції визначені при $S\in\Gamma\setminus \mathcal{X}$ та $Z\in\overline{\mathcal{U}}\setminus \mathcal{X}$.  Подамо функцію  \eqref{Om2-int}  у вигляді
		\begin{equation}
			\label{Omeg2-trans}
			\Omega_{2}^{\widetilde{\varphi}}\left(S,Z\right)=\frac{\widetilde{\varphi}\left(S\right)}{2}\,\frac{\left( d_{2}\left(S,Z\right)\bigl(\sigma_{1}'\left(S\right)-d_{1}\left(S,Z\right)\bigr)-
				d_{1}\left(S,Z\right)\bigl(\sigma_{2}'\left(S\right)-d_{2}\left(S,Z\right)\bigr) \right)}{\left(d\left(S,Z\right)\right)^2}\,.
		\end{equation}

		\vskip 2mm
		
\begin{lemma}
			\label{BikhIntToInequaliS}
			{\it Нехай виконуються умови\emph{:}
				\begin{itemize}
					\item[{ a)}] контурна  похідна  $\sigma'$  неперервна та не обертається в   нуль на  $\Gamma \setminus \mathcal{X}$, а  тілесна похідна  $\sigma'(Z)$ задовольняє оцінку \eqref{InEqTilConDeCorP}\emph{;}
					\item[{ b)}]
					при деякому $r_0\in \left(0, \mathcal{R}_{0}/4\right)$ для функції $\sigma'$
					виконуються оцінки \eqref{ots-SiDerLoc1}, \eqref{ots-SiDerLoc} і
					функція $d(S,Z)$ задовольняє оцінки
					\eqref{IneQLow} при $j=\overline{1,m}$ \emph{;}
					\item[{ c)}] функції $\varphi_{l} \colon
					\partial D_{\zeta}\setminus\Upsilon_{\zeta} \longrightarrow \mathbb{R}$, $l\in\{1,3\}$, з рівності \eqref{varphi13} неперервні на
					$\partial D_{\zeta}\setminus\Upsilon_{\zeta}$ і
					задовольняють оцінку
					\begin{equation}\label{phi-toCornPoin}
						|\varphi_l(\zeta)|  \le c\,\prod_{j=1}^m \|\zeta-\zeta_j \|^{-\alpha_j} \qquad \forall\,\zeta\in\partial D_{\zeta}\setminus\Upsilon_{\zeta},
					\end{equation}
					де    $\alpha_j <1$, якщо\,\, $\beta_j \in (-1,0)$, i
					$\alpha_j <1/\left(\beta_j+1\right)$, якщо\, \,
					$\beta_j \in [0,1]$, і стала\, $c$\, не залежить від $\zeta$.
				\end{itemize}
				Тоді:
				\begin{itemize}
					\item[{ 1)}]
					функції $\Omega_1^{\varphi_l}$, $\Omega_2^{\varphi_l}$ при  $l\in\{1,3\}$ задовольняють умови
					лем~\ref{contBikhIntCorPs}, \ref{BoDiffBikhIntCorPs}, накладені на функцію $\Omega$, в яких треба покласти 		\[\beta:=\max\left\{\max_{\beta_j \in \left[0,1\right]}\left\{\left(\beta_j + 1\right)\alpha_j\right\}, \,
					\max_{\beta_j \in \left(-1,0 \right)}\left\{\left(\beta_j + 1\right)\alpha_j - \beta_j\right\}
					\right\},\] 
					\[\beta_0:= \max \left\{0, -\beta_1, -\beta_2, \dots, -\beta_m\right\}\]
					і $\omega_0=\omega_1=\omega_2=\omega$, де функція $\omega$  визначається рівністю \eqref{DefOm-tsetrZ};
					\item[{ 2)}] функції $\Omega_1^{\varphi_l}$, $\Omega_2^{\varphi_l}$ при  $l\in\{1,3\}$ задовольняють умови
					леми~\ref{IneqOmegMoDiffPract}, накладені на функцію $\Omega$, в яких треба покласти 
					$\gamma_j=(\beta_j+1)\alpha_j$ та
					$\gamma_j'=-\beta_j$, якщо $\beta_j \in (-1,0)$, і $\gamma_j'=0$, якщо $\beta_j \in[0,1]$.
			\end{itemize} }
	\end{lemma}
	
	{\bf \em Доведення.} Для функції $\Omega=\Omega_1^{\varphi_l}$ при  $l\in\{1,3\}$ доведемо спочатку оцінки \eqref{omeg-loc-bound} і \eqref{omeg-locBound-circle}, а також оцінки
	\eqref{ineqOmonDomCp} і \eqref{ineqOmonBoundCp}, в яких замість $r_0$ слід покласти $r_0/2$.
	
	Враховуючи оцінки \eqref{phi-toCornPoin}, \eqref{Lemm8MainOts} і  \eqref{d-oz-znyzu-zag}, при
	$Z \in \overline{\mathcal{U}}\setminus \mathcal{X}$, $S\in \Gamma\setminus \mathcal{X}$ отримуємо
	\[ \left|\Omega_1^{\varphi_l} \left(S,Z\right) \right| =\left|\varphi_l\left(\widetilde{\sigma}(S)\right)\right|\frac{\left| \sigma'(S) - d\left(S,Z\right)\right|}{\left|d(S,Z) \right|} \le  c\, \prod_{j=1}^{m} \left|\sigma\left(S\right) - \sigma\left(X_j\right) \right|^{-\alpha_j}
	\min\left\{\frac{\omega\left(\left|S-Z\right|\right)}{\omega\left(\mathcal{R}\left(Z\right)/2\right)}, 1 \right\}\times\]
	\begin{equation}
		\label{sumAjBj}
		\times \Biggl(1+\prod_{j=1}^{m} \bigg(\frac{|S-X_j|}{|S-X_j|+|Z-X_j|} \bigg)^{\beta_j} \Biggr)\,,
	\end{equation}
	де стала\, $c$\, не залежить від $S$ і $Z$.
	
	Внаслідок оцінки \eqref{d-oz-znyzu-zag} маємо співвідношення
	\begin{equation}
		\label{varphiNearXj}
		\left|\sigma\left(S\right) - \sigma\left(X_j\right)  \right|^{\alpha_j}=\left|d\left(S,X_j\right) \right|^{\alpha_j}\left|S-X_j \right|^{\alpha_j} \ge \left|S-X_j \right|^{\left(\beta_j + 1\right)\alpha_j},\quad j=\overline{1,m}.
	\end{equation}
	
	Оцінимо зверху останній множник у правій частині нерівності \eqref{sumAjBj}.
	
	У випадку $S\in\Gamma\setminus \Gamma_{r_0/2}\langle \mathcal{X}\rangle$ існує додатна стала $M$, яка, взагалі кажучи, залежить від $r_0$, але не залежить від $S$ і $Z$, така, що
	\begin{equation}
		\label{Drib-oz}
		1+\prod_{j=1}^{m} \bigg(\frac{|S-X_j|}{|S-X_j|+|Z-X_j|} \bigg)^{\beta_j}\le M.
	\end{equation}
	
	Нехай тепер $S \in \Gamma_{r_0/2}(X_j)\setminus\{X_j\}$ при деякому $X_j\in\mathcal{X}$. Якщо $0<|Z-X_j|\le |S-X_j|$, то виконується оцінка вигляду \eqref{Drib-oz}. Якщо ж $|Z-X_j|\ge |S-X_j|$, то виконується оцінка
	\begin{equation}
		\label{Drib-oz-2}
		1+\prod_{j=1}^{m} \bigg(\frac{|S-X_j|}{|S-X_j|+|Z-X_j|} \bigg)^{\beta_j}\le
		\left \{ \begin{array}{ll}
			c\,\frac{|S-X_j|^{\beta_j}}{|Z-X_j|^{\beta_j}} & \mbox{при}\,\,\, \beta_j\in(-1,0), \\[3mm]
			c & \mbox{при}\,\,\, \beta_j\in [0,1],
		\end{array} \right.
	\end{equation}
	де стала\, $c$\,  не залежить від $S$ і $Z$.
	
	З нерівностей \eqref{sumAjBj} -- \eqref{Drib-oz-2} випливає виконання для функції $\Omega=\Omega_1^{\varphi_l}$ при  $l\in\{1,3\}$ нерівностей
	\eqref{omeg-loc-bound} і \eqref{omeg-locBound-circle}, а також нерівностей
	\eqref{ineqOmonDomCp} і \eqref{ineqOmonBoundCp}, в яких замість $r_0$ слід покласти $r_0/2$.

	Тепер для функції $\Omega=\Omega_1^{\varphi_l}$ при  $l\in\{1,3\}$ доведемо оцінки \eqref{omeg-locBound-diff}, \eqref{omeg-locBound-diffBound} і \eqref{ineqDiffOmBound}.
	
	Враховуючи оцінки \eqref{phi-toCornPoin}, \eqref{InEqTilConDeCorP}, \eqref{d-oz-znyzu-zag} і \eqref{IneqdifDZ0}, при
	$Z_0\in \Gamma \setminus \mathcal{X}$,  $Z \in \Gamma_{\frac{1}{4}\mathcal{R}(Z_0)}(Z_0)$ і $S \in \Gamma \setminus \mathcal{X}: |S-Z_0| \ge 2|Z-Z_0|$ отримуємо
	\[\left|\Omega_{1}^{\varphi_l}\left(S,Z\right)-
	\Omega_{1}^{\varphi_l}\left(S,Z_0\right) \right| =
	\left|\varphi_l\left(\widetilde{\sigma}\left(S\right)\right)\right|\left|\sigma'\left(S\right)\right|\frac{\left|d\left(S,Z\right) - d\left(S,Z_0)\right)\right|}{\left| d\left(S,Z\right)\right| \left| d\left(S,Z_0\right)\right|}\le\]
	\[ \le c\,  \prod_{j=1}^{m} \left|\sigma\left(S\right) - \sigma\left(X_j\right) \right|^{-\alpha_j}\, \frac{|Z-Z_0|}{|S-Z_0|}\,
	\min\left\{\frac{\omega\left(\left|S-Z_0\right|\right)}{\omega\left(\mathcal{R}\left(Z_0\right)/2\right)}, 1 \right\}\,\times  \]
	\begin{equation}
		\label{diffOm1ToL2}
		\times\, \prod_{j=1}^{m} \bigg(\frac{|S-X_j|}{|S-X_j|+|Z-X_j|} \bigg)^{\beta_j}\,
		\Biggl(1+\prod_{j=1}^{m} \bigg(\frac{|Z-X_j|+|Z_0-X_j|}{|S-X_j|+|Z_0-X_j|} \bigg)^{\beta_j} \Biggr)\,,
	\end{equation}
	де стала\, $c$\,  не залежить від $S$, $Z$ і $Z_0$.
	
	Для добутку останніх двох множників у правій частині нерівності \eqref{diffOm1ToL2} виконується оцінка
	\begin{multline}
		\label{Dob-drob-oz-2}
		\prod_{j=1}^{m} \bigg(\frac{|S-X_j|}{|S-X_j|+|Z-X_j|} \bigg)^{\beta_j}\,
		\Biggl(1+\prod_{j=1}^{m} \bigg(\frac{|Z-X_j|+|Z_0-X_j|}{|S-X_j|+|Z_0-X_j|} \bigg)^{\beta_j} \Biggr)\le\\[2mm]
		\le c\,\prod_{j : \beta_j<0} \Bigl(\min\left\{|S-X_j|, |Z_0-X_j|\right\}\Bigr)^{\beta_j},
	\end{multline}
	де стала\, $c$\,  не залежить від $S$, $Z$ і $Z_0$.
	
	З нерівностей \eqref{diffOm1ToL2}, \eqref{varphiNearXj} і \eqref{Dob-drob-oz-2} випливає виконання нерівностей \eqref{omeg-locBound-diff}, \eqref{omeg-locBound-diffBound} і \eqref{ineqDiffOmBound} для функції $\Omega=\Omega_1^{\varphi_l}$ при  $l\in\{1,3\}$.
	
	Розглянемо тепер функцію $\Omega=\Omega_2^{\varphi_l}$ при  $l\in\{1,3\}$  і доведемо для неї оцінки \eqref{omeg-loc-bound} і \eqref{omeg-locBound-circle}, а також оцінки
	\eqref{ineqOmonDomCp} і \eqref{ineqOmonBoundCp}, в яких замість $r_0$ слід покласти $r_0/2$.
	
	Враховуючи рівність \eqref{Omeg2-trans}, а також оцінки \eqref{phi-toCornPoin}, \eqref{sigmKIneq}, \eqref{d-oz-znyzu-zag} і  \eqref{Lemm8MainOts}, при
	$Z \in \overline{\mathcal{U}}\setminus \mathcal{X}$, $S\in \Gamma\setminus \mathcal{X}$ отримуємо
	\[ \left|\Omega_2^{\varphi_l} \left(S,Z\right) \right| =
	\frac{\Big|\varphi_l\left(\widetilde{\sigma}\left(S\right)\right)\Big|}{2}\,\frac{\Big|d_{2}\left(S,Z\right)\bigl(\sigma_{1}'\left(S\right)-d_{1}\left(S,Z\right)\bigr)-
		d_{1}\left(S,Z\right)\bigl(\sigma_{2}'\left(S\right)-d_{2}\left(S,Z\right)\bigr)\Big|}{\Big|d\left(S,Z\right)\Big|^2}\le\]
	\[ \le  c\, \prod_{j=1}^{m} \left|\sigma\left(S\right) - \sigma\left(X_j\right) \right|^{-\alpha_j}
	\min\left\{\frac{\omega\left(\left|S-Z\right|\right)}{\omega\left(\mathcal{R}\left(Z\right)/2\right)}, 1 \right\}\, \Biggl(1+\prod_{j=1}^{m} \bigg(\frac{|S-X_j|}{|S-X_j|+|Z-X_j|} \bigg)^{\beta_j} \Biggr)\,,\]
	де стала\, $c$\, не залежить від $S$ і $Z$.
	
	Отже, для функції $\Omega_2^{\varphi_l}$ виконується оцінка вигляду \eqref{sumAjBj}.
	Тому з нерівностей \eqref{varphiNearXj} -- \eqref{Drib-oz-2} випливає виконання для функції $\Omega=\Omega_2^{\varphi_l}$ при  $l\in\{1,3\}$ нерівностей
	\eqref{omeg-loc-bound} і \eqref{omeg-locBound-circle}, а також нерівностей
	\eqref{ineqOmonDomCp} і \eqref{ineqOmonBoundCp}, в яких замість $r_0$ слід покласти $r_0/2$.

	Доведемо нарешті для функції $\Omega=\Omega_2^{\varphi_l}$ при  $l\in\{1,3\}$ оцінки \eqref{omeg-locBound-diff}, \eqref{omeg-locBound-diffBound} і \eqref{ineqDiffOmBound}.
	
	Враховуючи рівність \eqref{Om2-int},
	при
	$Z_0\in \Gamma \setminus \mathcal{X}$,  $Z \in \Gamma_{\frac{1}{4}\mathcal{R}(Z_0)}(Z_0)$ і $S \in \Gamma \setminus \mathcal{X}: |S-Z_0| \ge 2|Z-Z_0|$ отримуємо
	\[\Omega_{2}^{\varphi_l}(S,Z)-\Omega_{2}^{\varphi_l}(S,Z_0)=\]
	\[=\frac{\varphi_l\left(\widetilde{\sigma}\left(S\right)\right)}{2}\,
	\bigg(\sigma'\left(S\right)
	\left(\frac{d_{2}\left(S,Z\right)}{\left(d\left(S,Z\right)\right)^2} -
	\frac{d_{2}\left(S,Z_0\right)}{\left(d\left(S,Z_0\right)\right)^2} \right) + 
	\sigma_{2}'\left(S\right)
	\left(
	\frac{1}{d\left(S,Z_0\right)} -
	\frac{1}{d\left(S,Z\right)} \right) \biggr)=\]
	\[=\frac{\varphi_l\left(\widetilde{\sigma}\left(S\right)\right)}{2} \,
	\biggl( \sigma_{2}'\left(S\right)\,
	\frac{d\left(S,Z\right)-d\left(S,Z_0\right)}{d\left(S,Z\right) d\left(S,Z_0\right)}   + \]
	\[+ \sigma'(S)\left(\frac{d_{2}\left(S,Z\right)- d_{2}\left(S,Z_0\right)}{\left(d_{2}\left(S,Z_0\right)\right)^2}+ \frac{d_2 \left(S,Z\right)\left(d\left(S,Z_0\right)- d\left(S,Z\right) \right) \left(d\left(S,Z_0\right) +  d\left(S,Z\right) \right)}{\left(d \left(S,Z\right)\right)^2 \left(d \left(S,Z_0\right)\right)^2}\right) \biggr). \]
	
	Звідси, враховуючи оцінки \eqref{phi-toCornPoin}, \eqref{InEqTilConDeCorP}, \eqref{d-oz-znyzu-zag} і \eqref{IneqdifDZ0},  отримуємо
	\[\left|\Omega_{2}^{\varphi_l}\left(S,Z\right)-
	\Omega_{2}^{\varphi_l}\left(S,Z_0\right) \right| \le \]
	\[\le c\, \left|\varphi_l\left(\widetilde{\sigma}(S)\right)\right|\,|\sigma'(S)|\, \frac{\left|d(S,Z)-d(S,Z_0)\right|}{\left|d\left(S,Z_0\right) \right|} \left(\frac{1}{\left|d\left(S,Z\right) \right|} +
	\frac{1}{\left|d\left(S,Z_0\right) \right|} \right)\le\]

	\[ \le c\,  \prod_{j=1}^{m} \left|\sigma\left(S\right) - \sigma\left(X_j\right) \right|^{-\alpha_j}\, \frac{|Z-Z_0|}{|S-Z_0|}\,
	\min\left\{\frac{\omega\left(\left|S-Z_0\right|\right)}{\omega\left(\mathcal{R}\left(Z_0\right)/2\right)}, 1 \right\}\,\times  \]
	\[	\times\, \prod_{j=1}^{m} \bigg(\frac{|S-X_j|}{|S-X_j|+|Z-X_j|} \bigg)^{\beta_j}\,
	\Biggl(1+\prod_{j=1}^{m} \bigg(\frac{|Z-X_j|+|Z_0-X_j|}{|S-X_j|+|Z_0-X_j|} \bigg)^{\beta_j} \Biggr)\,,\]
	де стала\, $c$\,  не залежить від $S$, $Z$ і $Z_0$.

	Отже, для функції $\Omega_2^{\varphi_l}$ виконується оцінка вигляду \eqref{diffOm1ToL2}.
	Тому з нерівностей  \eqref{varphiNearXj} і \eqref{Dob-drob-oz-2} випливає виконання нерівностей \eqref{omeg-locBound-diff}, \eqref{omeg-locBound-diffBound} і \eqref{ineqDiffOmBound} для функції $\Omega=\Omega_2^{\varphi_l}$ при  $l\in\{1,3\}$.
	
	Лему доведено.

	\vskip 2mm

	\begin{lemma}
		\label{FormsU13BikhIntKolCp}
		{\it Нехай виконуються умови леми~\ref{BikhIntToInequaliS}. Тоді для  функції \eqref{BihInt13ToSumOmeg},
де $\zeta=\widetilde{\sigma}(Z)$, при кожному  $Z_0\in \Gamma\setminus \mathcal{X}$ виконуються рівності
			\begin{equation}
				\begin{split}
					\label{LimComp13BikhonGamCornP}
					\lim\limits_{Z\to Z_0, Z \in \mathcal{U}}\mathrm{U}_{1}\left[	\Phi(\zeta) \right] = \frac{1}{2}\, \widetilde{\varphi}_{1} \left(Z_0\right)+ \mathrm{Re\,}\left(\mathcal{I}\left[\Omega_{\ast}^{\widetilde{\varphi}}\right] \left(Z_0\right)\right)+C_{\widetilde{\varphi}_{1}},\\
					\lim\limits_{Z\to Z_0, Z \in \mathcal{U}}\mathrm{U}_{3}\left[	\Phi(\zeta) \right] = \frac{1}{2}\, \widetilde{\varphi}_{3} \left(Z_0\right)+ \mathrm{Re\,}\left(\mathcal{I}\left[\Omega_{\ast\ast}^{\widetilde{\varphi}}\right] \left(Z_0\right)\right)+C_{\widetilde{\varphi}_{3}},
				\end{split}
			\end{equation}
де
	\begin{equation}
		\label{warpPhiK}
		\widetilde{\varphi}_{l}\left(Z_0\right) :=\varphi_{l} \left(\widetilde{\sigma}\left(Z_0\right)\right)\,, \qquad l\in\{1,3\},
	\end{equation}
	\begin{equation}
		\label{OmegAst1}
		\Omega_{\ast}^{\widetilde{\varphi}}\left(S,Z\right):= \Omega_1^{\widetilde{\varphi}_1}\left(S,Z\right)+2i
		\Omega_2^{\widetilde{\varphi}_1}\left(S,Z\right)-2\Omega_2^{\widetilde{\varphi}_3}\left(S,Z\right),
	\end{equation}
	\begin{equation}
		\label{OmegAst2}
		\Omega_{\ast \ast }^{\widetilde{\varphi}}\left(S,Z\right):=  \Omega_1^{\widetilde{\varphi}_3}\left(S,Z\right)- 2i
		\Omega_2^{\widetilde{\varphi}_3}\left(S,Z\right)-2\Omega_2^{\widetilde{\varphi}_1}\left(S,Z\right),
	\end{equation}
	\begin{equation}
		\label{intIomegVarph}
		\mathcal{I}\left[\Omega^{\widetilde{\varphi}}\right] \left(Z\right):= \frac{1} {2 \pi i } \int\limits_{\Gamma}
		\frac{\Omega^{\widetilde{\varphi}}\left(S,Z\right)}{S-Z} \, dS \qquad \forall\,Z \in \overline{\mathcal{U}}\setminus \mathcal{X},
	\end{equation}
\begin{equation}
		\label{int-const}
C_{\widetilde{\varphi}_l}:=\frac{1}{4 \pi i}\int\limits_{\Gamma} \frac{\widetilde{\varphi}_{l}\left(S\right)}{S}\, dS\,, \qquad l\in\{1,3\}.
\end{equation}
				}
	\end{lemma}
	
	{\bf \em Доведення.}
	З урахуванням співвідношень \eqref{Om1-int} -- \eqref{Omeg2-trans} і \eqref{warpPhiK} -- \eqref{intIomegVarph} та правил множення \eqref{tabl-umn} в бігармоніній алгебрі перепишемо рівність  \eqref{BihInt13ToSumOmeg} у вигляді
	\begin{multline*}
			\Phi(\zeta) = {e_1}\biggl(\mathcal{I}\left[\Omega_{\ast}^{\widetilde{\varphi}}\right] \left(Z\right) +
			 \frac{1}{2} \mathcal{S}\left[\widetilde{\varphi}_{1}\right]\left(Z\right)+C_{\widetilde{\varphi}_1}\biggr)+\\ +
			e_2 \biggl(\mathcal{I}\left[\Omega_{\ast \ast}^{\widetilde{\varphi}}\right] \left(Z\right) + \frac{1}{2} \mathcal{S}\left[\widetilde{\varphi}_{3}\right]\left(Z\right)+C_{\widetilde{\varphi}_3}\biggr)
			\qquad \forall\, Z \in \mathcal{U},
			\end{multline*}
	де
	\[\mathcal{S}\left[\widetilde{\varphi}_{l}\right]\left(Z\right):= \frac{1} {2 \pi i } \int\limits_{\Gamma} \frac{\widetilde{\varphi}_{l}\left(S\right)}{S}\frac{S+Z}{S-Z} \, dS
	\qquad  \forall\, Z \in \mathcal{U},\quad l\in\{1,3\},\,\,\, \mbox{---}\]
комплексний інтеграл Шварца.	
	
Оскільки сталі $C_{\widetilde{\varphi}_1}$ і $C_{\widetilde{\varphi}_3}$ --- дійсні, то
	\begin{equation}
		\begin{split}
			\label{Comp13BikhonGamCornP}
			\mathrm{U}_{1}\left[	\Phi(\zeta) \right] = \mathrm{Re\,}\left( \mathcal{I}\left[\Omega_{\ast}^{\widetilde{\varphi}}\right] \left(Z\right)+ \frac{1}{2}\,\mathcal{S}\left[\widetilde{\varphi}_{1}\right]\left(Z\right)\right)+C_{\widetilde{\varphi}_1},\\
			\mathrm{U}_{3}\left[	\Phi(\zeta) \right] = \mathrm{Re\,}\left( \mathcal{I}\left[\Omega_{\ast \ast}^{\widetilde{\varphi}}\right] \left(Z\right)+ \frac{1}{2}\,\mathcal{S}\left[\widetilde{\varphi}_{3}\right]\left(Z\right)\right)+C_{\widetilde{\varphi}_3}.
		\end{split}
	\end{equation}
	
Нехай $Z_0\in \Gamma\setminus \mathcal{X}$.	За умов леми~\ref{BikhIntToInequaliS}, враховуючи лему \ref{contBikhIntCorPs}, маємо рівності
\[\lim\limits_{Z\to Z_0, Z \in \mathcal{U}} \mathcal{I}\left[\Omega_{\ast}^{\widetilde{\varphi}}\right] \left(Z\right)= \mathcal{I}\left[\Omega_{\ast}^{\widetilde{\varphi}}\right] \left(Z_0\right), \qquad
	\lim\limits_{Z\to Z_0, Z \in \mathcal{U}} \mathcal{I}\left[\Omega_{\ast\ast}^{\widetilde{\varphi}}\right] \left(Z\right)= \mathcal{I}\left[\Omega_{\ast\ast}^{\widetilde{\varphi}}\right] \left(Z_0\right). \]
Крім того, для комплексного інтеграла Шварца виконуються рівності
\[\lim\limits_{Z\to Z_0, Z \in \mathcal{U}} \mathrm{Re\,}\, \mathcal{S}\left[\widetilde{\varphi}_{l}\right](Z)=\widetilde{\varphi}_{l}(Z_0), \qquad l\in\{1,3\}. \]
Тому, виконуючи граничний перехід у рівностях \eqref{Comp13BikhonGamCornP} при $Z\to Z_0$, $Z\in\mathcal{U}$, отримуємо рівності \eqref{LimComp13BikhonGamCornP}. Лему доведено.
	
	\vspace*{2mm}

При відшуканні розв'язків  (1-3)-задачі у класі функцій, що подаються у вигляді \eqref{cl-sol-cone},
припускаємо, що функції $\varphi_1$ і $\varphi_3$ з рівності \eqref{varphi13} неперервні на
$\partial D_{\zeta}\setminus\Upsilon_{\zeta}$ і задовольняють оцінку (\ref{phi-toCornPoin}).
З нерівності (\ref{varphiNearXj}) випливає, що відповідні їм функції (\ref{warpPhiK}) задовольняють нерівність
\begin{equation}\label{ozin-kl-K}
|\widetilde{\varphi}_{l}(S)| \le c \prod_{j=1}^{m} |S-X_j|^{-\gamma_j}\qquad \forall\, S \in \Gamma\setminus  \mathcal{X},
\end{equation}
де $\gamma_j=(\beta_j+1)\alpha_j<1$ і стала\, $c$\, не залежить від $S$.

Через $\mathcal{K}\{\gamma_j\}_{j=1}^m$ позначимо клас функцій $\widetilde{u} \colon \Gamma\setminus \mathcal{X} \longrightarrow \mathbb{R}$, що неперервні на $\Gamma\setminus \mathcal{X}$ і задовольняють оцінку (\ref{ozin-kl-K}), у якій усі $\gamma_j<1$ і стала\, $c$\, не залежить від $S$.

	Введемо  в розгляд функції, 
асоційовані з крайовими умовами   \eqref{Pr13week} кусково-неперервної (1-3)-задачі:
	\[\widetilde{u}_{l}\left(Z_0\right):= u_{l}\left(\widetilde{\sigma}\left(Z_0\right)\right)\qquad \forall\, Z_0 \in \Gamma
	\setminus \mathcal{X}, \quad l\in\{1,3\}.\]

Основним результатом роботи є наступна

\vskip 2mm

\begin{theorem}
			\label{Th-main}
{\it Нехай межа області $D_{\zeta}$ така, що виконуються умови a) і b) леми~\ref{BikhIntToInequaliS}, і нехай функції $u_{j} \colon \partial D_{\zeta} \setminus \Upsilon_{\zeta} \longrightarrow \mathbb{R}$, $j\in \{1,3\}$, --- неперервні на $\partial D_{\zeta}\setminus\Upsilon_{\zeta}$ і задовольняють оцінку вигляду \eqref{phi-toCornPoin}.
Тоді розв'язання кусково-неперервної \em (1-3)\em-задачі у класі функцій, що подаються у вигляді \eqref{cl-sol-cone}, при цьому
функції $\varphi_1$ і $\varphi_3$ з рівності \eqref{varphi13} неперервні на
$\partial D_{\zeta}\setminus\Upsilon_{\zeta}$ і задовольняють оцінку \eqref{phi-toCornPoin}, зводиться до розв'язання у класі функцій $\mathcal{K}\{\gamma_j\}_{j=1}^m$ системи інтегральних рівнянь
	\begin{equation}
		\label{BVC13-circCorPo}
				\begin{array}{ll}
			\frac{1}{2}\,\widetilde{\varphi_1}\left(Z_0\right)+
			\mathrm{Re\,}\left(\mathcal{I}\left[\Omega_{1}^{\widetilde{\varphi}_1}\right] \left(Z_0\right)\right)- 2\, \mathrm{Im\,}\left(\mathcal{I}\left[\Omega_{2}^{\widetilde{\varphi}_1}\right] \left(Z_0\right)\right)-\\[2mm]
			\hspace*{45mm} -2\, \mathrm{Re\,}\left(\mathcal{I}\left[\Omega_{2}^{\widetilde{\varphi}_3}\right] \left(Z_0\right)\right)+C_{\widetilde{\varphi}_1}=
			\widetilde{u}_{1}\left(Z_0\right),\\[4mm]
			\frac{1}{2} \, \widetilde{\varphi_3}\left(Z_0\right)+
			 \mathrm{Re\,}\left(\mathcal{I}\left[\Omega_{1}^{\widetilde{\varphi}_3}\right] \left(Z_0\right)\right)+ 2 \, \mathrm{Im\,}\left(\mathcal{I}\left[\Omega_{2}^{\widetilde{\varphi}_3}\right] \left(Z_0\right)\right)-\\[2mm]
  \hspace*{45mm} -2\mathrm{Re\,}\left(\mathcal{I}\left[\Omega_{2}^{\widetilde{\varphi}_1}\right] \left(Z_0\right)\right)+C_{\widetilde{\varphi}_3}=\widetilde{u}_{3}\left(Z_0\right) \qquad \forall\, Z_0 \in \Gamma\setminus \mathcal{X},
		\end{array}
	\end{equation}
	де інтеграл $\mathcal{I}$ визначений рівністю \eqref{intIomegVarph},
функції $\widetilde{\varphi}_{l}$ і сталі $C_{\widetilde{\varphi}_l}$ при $l\in \{1,3\}$  визначені відповідно рівностями \eqref{warpPhiK} і
\eqref{int-const}. }
	\end{theorem}
	
	{\bf \em Доведення.} Враховуючи рівності \eqref{LimComp13BikhonGamCornP}, \eqref{OmegAst1} і \eqref{OmegAst2}, переписуємо крайові умови
\eqref{Pr13week} кусково-неперервної (1-3)-задачі у вигляді системи інтегральних рівнянь \eqref{BVC13-circCorPo} для відшукання функцій $\widetilde{\varphi}_{1}, \widetilde{\varphi}_{3}\in \mathcal{K}\{\gamma_j\}_{j=1}^m$, де задані функції $\widetilde{u}_{1}$ і $\widetilde{u}_{3}$ також належать класу $\mathcal{K}\{\gamma_j\}_{j=1}^m$.

З лем \ref{BoDiffBikhIntCorPs} і \ref{IneqOmegMoDiffPract} випливає, що клас функцій $\mathcal{K}\{\gamma_j\}_{j=1}^m$ є інваріантним відносно інтегральних операторів із лівих частин рівнянь \eqref{BVC13-circCorPo}. Отже,
 можна стверджувати, що  розв'язання кусково-неперервної (1-3)-задачі у класі функцій, що подаються у вигляді \eqref{cl-sol-cone} з
 наведеними в теоремі припущеннями про функції $\varphi_1$ і $\varphi_3$ з рівності \eqref{varphi13}, еквівалентне
 розв'язанню системи інтегральних рівнянь \eqref{BVC13-circCorPo} у класі $\mathcal{K}\{\gamma_j\}_{j=1}^m$. Теорему доведено.

 \vskip 2mm

{\bf Зауваження.}

{\bf 1.} Серед областей $D_{\zeta}$, що задовольняють умови теореми \ref{Th-main}, є області, межа яких в точках набору $\Upsilon_{\zeta}$ може утворювати нульові кути з вістрям, направленим всередину області. Такою, наприклад, є область, конгруентний "двійник"\/ якої $D_z$ у комплексній площині є образом одиничного круга $\mathcal{U}$ при конформному відображенні $\sigma(Z)=(Z+1)^2$.

{\bf 2.} Якщо в усіх інтегралах, що входять в рівняння системи \eqref{BVC13-circCorPo}, інтегрування по колу $\Gamma$ замінити інтегруванням по множині $\Gamma\setminus\Gamma_{\delta}\langle \mathcal{X}\rangle$ при достатньо малому $\delta>0$, то в силу леми \ref{BoDiffBikhIntCorPs} отримані інтеграли задаватимуть компакні оператори в банаховому просторі неперервних на $\Gamma\setminus\Gamma_{\delta}\langle \mathcal{X}\rangle$ функцій при стандартному визначенні норми такого простору. Отже, рівняння, отримані з рівнянь системи \eqref{BVC13-circCorPo} в результаті зазначеної заміни множини інтегрування, будуть фредгольмовими у банаховому просторі неперервних на $\Gamma\setminus\Gamma_{\delta}\langle \mathcal{X}\rangle$ функцій при усіх достатньо малих $\delta>0$. У той же час, питання про збіжність розв'язків вказаних рівнянь при $\delta\to 0+0$ у належних вагових просторах функцій класу $\mathcal{K}\{\gamma_j\}_{j=1}^m$ потребує спеціального дослідження.


\begin{thebibliography}{99}
		
		\bibitem{Mikhlin}
		С.\,Г.~Михлин, {\em Плоская задача теории упругости}, Труды сейсм. ин-та АН СССР, №~65 (1934), 83~c.
		
		\bibitem{Kondrat'ev-TMMO-67}
		В.\,А. Кондратьев, {\em Краевые задачи для эллиптических уравнений в областях с коническими или угловыми точками}, Тр. ММО, \textbf{16}, Издательство Московского университета, Мoсква,  209--292 (1967).
		
	
		\bibitem{Kufner-Sandig}
		A. Kufner,  A.-M. S\"{a}ndig,  {\em  Some applications of weighted Sobolev spaces} (TEUBNER-TEXTE
		zur Mathematik, \textbf{100}),
		BSB BG Teubner Verlagsgesellschaft, Leipzig (1987).
		
		\bibitem{MaNazPlam-ASumptTh-V1}
		V. Maz'ya,  S. Nazarov,  B. Plamenevskij, {\em   Asymptotic theory of elliptic boundary value problems in singularly perturbed domains. Vol.~1} (Operator Theory
		Advances and Applications, \textbf{111}), Springer Science \& Business Media (2000).
		
		\bibitem{Mush_upr}
		Н.\,И. Мусхелишвили, {\em   Некоторые основные задачи математической теории
			упругости}, Изд-во "Наука", Москва (1966).
		
		\bibitem{Lurie_engl}  Lurie A.I.
		\emph{ Theory of Elasticity. Engl. transl. by A. Belyaev}.
		Springer-Verlag: Berlin etc., 2005.
		
		
		\bibitem{Mikh_int_eq_Th-El}
		S.\,G. Mikhlin,  N.\,F. Morozov, M.\,V. Paukshto,  {\em The integral equations of the theory of elasticity} (TEUBNER-TEXTE
		zur Mathematik, \textbf{135}), Springer, Stuttgart – Stiint. etc. (1995).
		
		\bibitem{Mikh_kniga} Mikhlin S.G. \emph{ Integral Equations and their Applications to
			Certain Problems in Mechanics, Mathematical Physics and
			Technology}. Pergamon Press: New York, 1964.
		
	
		\bibitem{Kupragze63}
		В.\,Д. Купрадзе, {\em Методы потенциала в теории упругости}, Гос. изд-во физю-мат. лит-ры, Москва (1963).
		
		
		\bibitem{Lopatinskii}
		Я.\,Б. Лопатинский,   {\em Об одном способе приведения граничных задач для системы дифференциальных уравнений эллиптического типа к регулярным интегральным уравнениям}, Укр. мат. журн, \textbf{5}, № 2, 123-151 (1953).
		
	
		\bibitem{Maz'ya91Engl}
		В.\,Г. Мазья, {\em Граничные интегральные уравнения}, Анализ-4, Итоги науки и техн. (Сер. Соврем. пробл. мат. Фундам. направления), \textbf{27}, ВИНИТИ,  Мoсква,  131-–228 (1988).
		
		\bibitem{Magnaradze-38}
		Л.\,Г. Магнарадзе,  {\em Основные задачи плоской теории упругости для контуров с угловыми точками}, Труды Тбилисского математического института, \textbf{4},  Тбилиси, 43--76 (1938).
		
		\bibitem{Radon-46}
		И.\,О. Радон,   {\em О краевых задачах для логарифмического потенциала}, Успехи мат. наук, \textbf{1},
		№3-4(13-14), 96--124 (1946).
		
		\bibitem{Polozhyi}
		Г.\,Н. Положий, {\em  Решение некоторых задач плоской теории упругости для областей с угловыми точками},  Укр. мат.  журн., \textbf{1}, № 4, 16--41 (1949).
		
		\bibitem{G-Albinus-83PAN}
		G. Albinus, {\em Multiple layer potentials for the quadrant and their application to the Dirichlet problem in plane domains with a piecewise smooth boundary},   Banach Center Publications, \textbf{10}, №1,  Warsaw, 7--26 (1983).
		
		\bibitem{Meln-PLaksa-PBP-quadrante}
		Мельниченко И.П., Плакса С.А., {\em Редукция основной бигармонической задачи для квадранта к несингулярным интегральным уравнениям},  Укр. мат.  журн., \textbf{47}, №~6, 775--784  (1995).
		
	
		\bibitem{GrPl-Angle-umz}
		Грищук~С.\,В.,   Плакса~С.\,А.,  {\em Бігармонічна задача для кута і моногенні функції},
		Укр. мат. журн.  \textbf{74}, №~11, C.~1478--1491 (2022);\\
		{\bf Engl. translation:}
		S. V. Gryshchuk and S. A. Plaksa,
		{\em Biharmonic problem for an angle and monogenic functions},
		Ukr. Math. J.,  \textbf{74}, No.~11
		Springer Science+Business Media, LLC, pp.~1686--1700 (2023).
		
		
		
		\bibitem{IJPAM_13}
		S.\,V. Gryshchuk, S.\,A. Plaksa, {\em   Schwartz-type integrals in a biharmonic
			plane},  Intern. J. of Pure and Appl. Math.,  \textbf{83}, № 1,
		193--211 (2013).
		
		
		\bibitem{mon-f-bih-BVP-MMAS}
		S.\,V.~Gryshchuk, S.\,A. Plaksa, {\em Monogenic functions in the
			biharmonic boundary value problem}, Math. Methods Appl. Sci.,
		\textbf{39}, No.~11, 2939--2952 (2016).
		
		\bibitem{Gr-1Fr17}
		С.\,В. Грищук, {\em  Одновимірність ядра системи інтегральних рівнянь Фредгольма для однорідної бігармонічної задачі}, Зб. праць Ін-ту математики НАН України, \textbf{14}, №1, 128--139 (2017).
		
		\bibitem{BeghGrPl19}
		Gryshchuk~S.\,V., Plaksa~S.\,A., {\em   Schwartz-Type Boundary Value Problems for Monogenic Functions in a Biharmonic Algebra}. In:  S. Rogosin, A. \c{C}elebi  (eds.),  {Analysis as a Life, Dedicated to Heinrich Begehr on the Occasion of his 80th Birthday (Trends in Mathematics)},  Birkha\"{u}ser,  193--211 (2019).
		
		
		\bibitem{GrPlHBVMbF-Flaut21}
		Gryshchuk~S.\,V., Plaksa~S.\,A., {\em A Hypercomplex Method for Solving Boundary Value Problems for Biharmonic Functions}. In: \v{S} Ho\v{s}kov\'a-Mayerov\'a,  C. Flaut, F. Maturo  (eds), Algorithms as a Basis of Modern Applied Mathematics (Studies in Fuzziness and Soft Computing), \textbf{404}, Springer, Cham, 231--255 (2021).
		
		
		
		\bibitem{SrPlCanD-UMB21}
		S.\,V.~Gryshchuk, S.\,A. Plaksa, {\em Schwartz-type boundary value problems for
			canonical domains in a biharmonic plane}, Ukrainian Mathematical Bulletin, \textbf{18}, №3, 338--358 (2021).
		
		
		\bibitem{KM-BFf}
		В. Ф. Ковалев, И. П. Мельниченко, {\em  Бигармонические функции на бигармонической плоскости}, Доп. АН УРСР.
		Сер. А, № 8, 25--27 (1981).
		
		\bibitem{Mel86}
		И.\,П.~Мельниченко, {\em    Бигармонические базисы в алгебрах второго
			ранга}, Укр. мат. журн., {\bf 38}, №2, 252--254 (1986).
		
		
		
		\bibitem{Sodbero}
		L. Sobrero,   {\em Nuovo metodo per lo studio dei problemi di
			elasticit\`{a}, con applicazione al problema della piastra forata},
		Ricerche di Ingegneria,  \textbf{13}, №2, 255--264 (1934).
		
		\bibitem{Douglis-53}
		A.  Douglis, {\em A function-theoretic approach to elliptic systems of equations in
			two variables},  Communications on Pure and Applied Mathematics,
		\textbf{6}, №2, 259--289 (1953).
		
		\bibitem{GrPl_umz-09}
		С.\,В. Грищук, С.\,А. Плакса,
		{\em Моногенные функции в бигармонической
			алгебре}, Укр. мат. журн., {\textbf{61}},  №12, 1587--1596  (2009).
		
		
		\bibitem{dopovidi_09}
		Грищук С.В., Плакса С.А., {\em Моногенные функции в бигармонической плоскости}, Доповіді НАН України, Сер. Мат. прир. техн. науки, № 12., C.~13--20 (2009).
		
		
		
		\bibitem{Conrem_13}
		S.\,V. Gryshchuk, S.\,A. Plaksa, {\em Basic Properties of Monogenic Functions
			in a Biharmonic Plane}. In:  M.\,L.  Agranovskii,  M. Ben-Artzi,
		G. Galloway, L. Karp,  V. Maz'ya (eds.), ``Complex Analysis and
		Dynamical Systems V''. Contemporary Mathematics,  \textbf{591},
		Amer. Math. Soc., Providence, RI, 127--134 (2013).
		
		
		
		\bibitem{Kovalov}
		В.\,Ф. Ковалев, {\em  Бигармоническая задача Шварца},  Киев, (1986), 19 с. (Препр. / НАН Украины.
		Ин-т математики; 86.16).
		
		\bibitem{Walsh-Sewell-1940}
		J. L. Walsh,  W. E. Sewell,
		{\em Sufficient conditions for various degrees of approximation by polynomials}, Duke Math. J., {\textbf{6}}, 658–-705 (1940).
		
		
		
		\bibitem{tamrazov}
		П. М. Тамразов, Гладкости и полиномиальные приближения, Наукова думка, Киев, 1975. 
		
		
		\bibitem{Dz-Shev-6-3}
		Dzyadyk~V.\,K., Shevchuk~I.\,A., {\em Theory of uniform approximation of functions by polynomials}. De Gruyter, Berlin (2008).
		
	\end{thebibliography}
\end{document}